\documentclass[11pt,twoside]{amsart} 
\usepackage{amssymb}
\usepackage{amsmath}
\usepackage[OT2,T2A]{fontenc}
\advance \topmargin by -\headheight
\advance \topmargin by -\headsep
\evensidemargin \oddsidemargin
\date{}
\newtheorem{theorem}{Theorem}[section]
\newtheorem{conjecture}[theorem]{Conjecture}
\newtheorem{corollary}[theorem]{Corollary}
\newtheorem{lemm}{Lemma}

\newtheorem{lemma}[theorem]{Lemma}
\newtheorem{prop}{Proposition}

\newtheorem{proposition}[theorem]{Proposition}
\title{Smooth representations and sheaves} 
\author{U.Jannsen, M.Rovinsky} 
\address{NWF I-Mathematik, Universit\"{a}t Regensburg, 
Universit\"{a}tsstra{\ss}e 31, 93053 Regensburg}
\address{Independent University of Moscow, 119002 Moscow 
B.Vlasievsky Per. 11 \& Institute for Information 
Transmission Problems of Russian Academy of Sciences}
\dedicatory{To Pierre Deligne on the occasion of his 65th birthday}
\begin{document} 
\begin{abstract} The paper is concerned with `geometrization' of smooth 
(i.e. with open stabilizers) representations of the automorphism group of 
universal domains, and with the properties of \ `geometric' representations 
of such groups. As an application, we calculate the cohomology groups of 
several classes of smooth representations of the automorphism group of 
an algebraically closed extension of infinite transcendence 
degree of an algebraically closed field. 
\end{abstract}

\maketitle 

Our motivation comes from an attempt to understand interplay between geometry 
($K$-groups) and topology (cohomology) of algebraic varieties, 
expressed in the motivic conjectures, 
through representations of certain `big' totally disconnected topological 
groups, cf. \cite{rms}. However, in this paper we reverse the point of view 
and study `geometrization' of such representations. 

\subsection{Notations and terminology} Let $k$ be a characteristic zero 
field, $F|k$ be an algebraically closed extension of countable transcendence 
degree, and $G=G_{F|k}$ 
be the group of field automorphisms of $F$ leaving $k$ fixed. Let $B$ be the 
collection of the pointwise stabilizers of the finite subsets in $F$, i.e., 
the subgroups of type $G_{F|L}$ for subfields $L$ of $F|k$ of finite type. 

We consider $G$ as a totally disconnected topological group with base of open 
subgroups $B$. These notions are recalled in Appendix \ref{general-topology}. 

For a totally disconnected group $H$ an $H$-set is called {\sl smooth} if 
the stabilizers of all its elements are open. 

Let $E$ be a field endowed with a smooth $H$-action 
(by field automorphisms). 

Denote by ${\mathcal S}m_H(E)$ the category of 
smooth $E$-semilinear representations of $H$, i.e., $E$-vector spaces $V$ 
endowed with a smooth additive action of $H$ such that $h(ev)=he\cdot hv$ 
for any $e\in E$, $h\in H$, $v\in V$. 
We set ${\mathcal S}m_H:={\mathcal S}m_H({\mathbb Q})$. Here the $H$-action 
on ${\mathbb Q}$ is trivial, which is the only possibility in the case $H=G$. 

The full subcategory ${\mathcal I}_G$ of ${\mathcal S}m_G$, whose objects 
are `homotopy invariant', is defined on p.\pageref{homotopy-invariant}. 

For a variety $X$ over a field we denote by $CH^q(X)$ the (Chow) group 
of cycles of codimension $q$ modulo rational equivalence, i.e., modulo 
divisors of rational functions on subvarieties of codimension $q-1$. 

For an abelian group $A$ we set $A_{{\mathbb Q}}:=A\otimes{\mathbb Q}$. 

For a unitary commutative ring $A$ and $A$-algebras $A_1$ and $A_2$ denote by 
$\{A_1\stackrel{/A}{\longrightarrow}A_2\}$ the set of unitary $A$-homomorphisms 
$A_1\stackrel{/A}{\longrightarrow}A_2$. (If $A_1$ is a field then 
the homomorphisms are injectives, so the notation will be 
$\{A_1\stackrel{/A}{\hookrightarrow}A_2\}$.) If $S$ is a set then $A[S]$ 
denotes the free $A$-module with basis $S$. 

An algebraic closure of a field $L$ is denoted by $\overline{L}$, if $L$ is 
a subfield in $F$ then $\overline{L}$ denotes the algebraic closure in $F$. 

\subsection{Context and content} We are interested in interplay between 
algebraic geometry and representation theory of $G$. Various geometric 
categories are linked to the representations of $G$ by fully faithful functors. 
An example of such full embedding is given in Appendix \ref{ind-comp}, 
more examples are in 
\cite[Theorem 1.1, Proposition 4.3]{repr}. It is explained 
in \cite[Corollary 7.9]{pgl} that it is especially important for 
geometric applications to describe the {\sl irreducible smooth} 
representations of $G$ (or at least `homotopy invariant' 
ones among them, cf. p.\pageref{homotopy-invariant} below). 

It is well-known that for any totally disconnected group $H$ the 
smooth (i.e. with open stabilizers) $H$-sets and their $H$-equivariant 
maps form a topos, cf. Appendix \ref{general-constr-sheaves}. 

In this paper, we introduce two sites, called the dominant and the Hartogs 
topologies, respectively. They are quite different from 
that described in \S\ref{general-constr-sheaves} in the case of $H=G$. 

In a way similar to Galois--Grothendieck correspondence, Lemma 
\ref{ekviv-st} describes the smooth $G$-sets as the sheaves in the 
dominant topology. Then we apply the standard sheaf-theoretic techniques 
to the study of the cohomology of smooth $G$-modules. Roughly, one of our 
results (Corollary \ref{2-3-acikl}) states the acyclicity of the smooth 
$G$-module corresponding to the dominant sheafification of a presheaf 
on the category of smooth $k$-varieties. 

We deduce in Corollary \ref{cohom-dim} that for $E={\mathbb Q}$ or 
$E=F$ (with the standard $G$-action) the category ${\mathcal S}m_G(E)$ 
is of infinite cohomological dimension, which is quite expectable. 

However, there are `too many' smooth representations of $G$, cf. Appendix 
\ref{count-cyclic-irr}. For this reason, instead of studying {\sl all} 
smooth representations of $G$, one can try to restrict the category 
of considered representations. More particularly, one can look for the 
irreducible $F$-semilinear representations of $G$ containing `homotopy 
invariant' representations, expecting that there are `very few' of them. 

To explore this option we introduce in \S\ref{geom-puchki} a functor 
from the category of smooth representations of $G$ to the category of 
sheaves in the Hartogs topology. From the point of vue of this functor, 
the `homotopy invariant' representations `look as local systems' 
(Proposition \ref{loc-i-glob-invar}), and the `interesting' 
$F$-semilinear representations of $G$ are `globally generated' 
(Remark 1 after Lemma \ref{basic-glob}). 

Homological machinery could be applied to several classes of questions. 
We mention just four examples. 
\begin{enumerate} \item The acyclicity of some classes of representations 
is often crucial, e.g. for the knowledge of the structure of some 
co-induced representations, cf. \S\ref{alg-and-subgroups}. 
\item
One knows several cases when the ${\rm Hom}$-groups in the triangulated 
category of mixed motives turn out to be related (dual) to corresponding 
${\rm Ext}$-groups in the category of `homotopy invariant' 
representations of $G$, cf. \cite{repr}. So one of the problems 
is to find other ${\rm Ext}$-groups in the category of `homotopy 
invariant' (and more generally, smooth) representations of $G$ 
and compare them with the conjectural values of corresponding 
${\rm Ext}$-groups in the conjectural category of mixed motives. 
\item
There are some conjectures coming from geometry: on the semi-simplicity 
of the graded quotients of the level filtration $N_{\bullet}$, cf. 
\cite[Conjecture 6.9]{repr}, on the splitting of the filtration 
${\mathcal F}^{\bullet}$, ibid, etc. 
\item
It is desirable to enlarge the category of `homotopy invariant' 
representations and relate it to the category of {\sl effective} mixed 
motives, so that in particular, this bigger category would contain objects 
like $F^{\times}$ and ${\mathbb Q}$ was still a projective object. 
\end{enumerate}

\subsection{From modules over Hecke algebras to modules over 
subquotient groups of $G$} \label{alg-and-subgroups}
In the standard type of semi-simplicity or irreducibility criteria of 
$G$-modules of Appendix \ref{Hecke-algebras}, it would be natural to 
replace the semi-simplicity or the irreducibility conditions for the 
modules over the Hecke algebras by the corresponding conditions on 
representations of the {\sl groups} $G_{F'|k}$ for the algebraically 
closed extensions $F'|k$ in $F$ of finite transcendence degree. This is 
how the semi-simplicity of the $G$-module $B^q(X_F)$ is established in 
\cite[Proposition 3.8]{repr}. (Here $B^q(X_F)$ is 
the Chow group $CH^q(X_F)_{{\mathbb Q}}$ modulo the 
`numerical equivalence over $k$'.) More precisely, this is based 
on the following 
\begin{lemm} Let $W\in{\mathcal S}m_G$. If 
${\rm Hom}_G(Z^{\dim X}(k(X)\otimes_kF),W)={\rm Hom}_G(CH_0(X_F),W)$ for 
any smooth proper variety $X$ over $k$ then $W$ is semi-simple if and only 
if the $G_{F'|k}$-module $W^{G_{F|F'}}$ is semi-simple for any algebraically 
closed $F'$ of finite transcendence degree over $k$. \end{lemm} 
{\it Proof.} Clearly, $W^{G_{F|L'}}=W^{G_{F|L}}$ for any purely 
transcendental extension $L'|L$ in $F|k$, i.e., $W\in{\mathcal I}_G$, 
cf. p.\pageref{homotopy-invariant} and \cite{repr}. 

By the semi-simplicity criterion of Appendix \ref{Hecke-algebras}, Remark 1, 
p.\pageref{semi-simplicity-general}, the $G$-module $W$ is semi-simple 
if and only if for any $L'|L$ as above with $L$ of finite type over $k$ 
and $F$ algebraic over $L'$ the ${\mathcal H}_{G_{F|L'}}$-module $W^{G_{F|L'}}$ 
is semi-simple. Here ${\mathcal H}_{G_{F|L'}}:=
h_{L'}{\mathbb D}_{{\mathbb Q}}(G)h_{L'}\supseteq\langle h_{L'}\sigma 
h_{L'}~|~\sigma\in G\rangle_{{\mathbb Q}}$ is the Hecke algebra and $h_{L'}$ 
is the Haar measure on $G_{F|L'}$. 

As $W$ is a quotient of a direct sum of objects of type $CH_0(X_F)_{{\mathbb Q}}$, 
the action $${\mathbb D}_{{\mathbb Q}}(G)\otimes W^{G_{F|L}}
{\longrightarrow\hspace{-3mm}\to}{\mathbb Q}[G/G_{F|L}]\otimes 
W^{G_{F|L}}\longrightarrow W$$ factors through $CH_0(Y_F)_{{\mathbb Q}}
\otimes W^{G_{F|L}}\longrightarrow W$, where $Y$ is a smooth 
proper model of the extension $L|k$, cf. \cite[Proposition 3.6]{repr}. 

Thus, the action of ${\mathcal H}_{G_{F|L'}}$ factors through 
$CH_0(Y_{k(Y)})_{{\mathbb Q}}=h_{L'}CH_0(Y_F)_{{\mathbb Q}}$. 

In other words, the action of the Hecke algebra ${\mathcal H}_{G_{F|L'}}(G)$ 
on $W^{G_{F|L}}$ is determined by the action of the Hecke algebra 
${\mathcal H}_{G_{\overline{L}|L}}(G_{\overline{L}|k})$, 
so the semi-simplicity of the ${\mathcal H}_{G_{F|L'}}(G)$-module 
$W^{G_{F|L}}$ is equivalent to its semi-simplicity as a 
${\mathcal H}_{G_{\overline{L}|L}}(G_{\overline{L}|k})$-module. \qed 

\vspace{4mm}

The following lemma reduces the general 
semi-simplicity problem to an acyclicity question. 
\begin{lemm} Let ${\mathcal H}$ be a subcategory in ${\mathcal S}m_G$ 
closed under taking subobjects and $F'|k$ be an algebraically closed 
extension in $F$ of a finite transcendence degree $q$. The following 
conditions on the subcategory ${\mathcal H}$ and on $F'$ are equivalent. 
\begin{enumerate} \item \label{peresech-s-invar} 
For any $W\in{\mathcal H}$ any $G_{F'|k}$-submodule 
$U\subseteq W^{G_{F|F'}}$ coincides with the $G_{F'|k}$-submodule 
of $G_{F|F'}$-invariants in the $G$-submodule generated by $U$: 
$U=\langle U\rangle_G^{G_{F|F'}}$. 
\item \label{sur-inv-usl} 
For any $W\in{\mathcal H}$, any finitely generated field extension $L|k$ 
of transcendence degree $q$ and any surjection ${\mathbb Q}
[\{L\stackrel{/k}{\hookrightarrow}F\}]^N\longrightarrow W$ 
in ${\mathcal S}m_G$ induces a surjection of $G_{F|F'}$-invariants 
${\mathbb Q}[\{L\stackrel{/k}{\hookrightarrow}F'\}]^N\longrightarrow 
W^{G_{F|F'}}$ in ${\mathcal S}m_{G_{F'|k}}$, where $F'$ is algebraic over $L$. 
\item \label{H-1-usl} For any finitely generated field extension $L|k$ 
of transcendence degree $q$ and for any $Q\subseteq{\mathbb Q}[\{L
\stackrel{/k}{\hookrightarrow}F\}]^N$ such that the quotient belongs 
to ${\mathcal H}$, one has $H^1_{{\mathcal S}m_G}(G_{F|F'},Q)=0$. 
\end{enumerate} \end{lemm} 
{\it Proof.} (\ref{sur-inv-usl}) and (\ref{H-1-usl}) are equivalent 
by Corollary \ref{cohom-dim}, since the restriction of ${\mathbb Q}
[\{L\stackrel{/k}{\hookrightarrow}F\}]$ to $G_{F|F'}$ is canonically 
isomorphic to ${\mathbb Q}[\{L\otimes_kF'\stackrel{/F'}
{\longrightarrow}F\}]=\bigoplus_{x\in{\rm Spec}(L\otimes_kF')}
{\mathbb Q}[\{F'(x)\stackrel{/F'}{\hookrightarrow}F\}]$. 

(\ref{sur-inv-usl})$\Rightarrow$(\ref{peresech-s-invar}). 
If $U\neq\langle U\rangle_G^{G_{F|F'}}$ then there are elements 
$u_1,\dots,u_N\in U$, $a_{ij}\in{\mathbb Q}$, and $g_{ij}\in G$ such that 
$\sum_{1\le i\le N,~j}a_{ij}g_{ij}u_i\in W^{G_{F|F'}}
\smallsetminus U$. Replace $W$ by the $G$-span of $u_1,\dots,u_N$, and 
$U$ by the $G_{F'|k}$-span of $u_1,\dots,u_N$. Then there is a surjection 
${\mathbb Q}[\{L\stackrel{/k}{\hookrightarrow}F\}]^N
\longrightarrow\hspace{-3mm}\to W$, $e_i\mapsto u_i$, where $L|k$ is an 
extension in $F'$ of finite type such that $G_{F|L}$ is contained in the 
common stabilizer of the elements $u_1,\dots,u_N$. Then (\ref{sur-inv-usl}) 
implies the surjectivity of ${\mathbb Q}[\{L\stackrel{/k}{\hookrightarrow}
F'\}]^N\longrightarrow\hspace{-3mm}\to W^{G_{F|F'}}$. As it factors through 
$U$, we get $U=W^{G_{F|F'}}$, as required. 

(\ref{peresech-s-invar})$\Rightarrow$(\ref{sur-inv-usl}). Let $U$ be the 
image of ${\mathbb Q}[\{L\stackrel{/k}{\hookrightarrow}F'\}]^N$, i.e., the 
$G_{F'|k}$-span of $u_1,\dots,u_N$. Then the $G$-span of $u_1,\dots,u_N$ 
coincides with $W$. It follows from (\ref{peresech-s-invar}) that 
$U=W^{G_{F|F'}}$, i.e., the homomorphism ${\mathbb Q}[\{L\stackrel{/k}
{\hookrightarrow}F'\}]^N\longrightarrow\hspace{-3mm}\to W^{G_{F|F'}}$ 
is surjective. \qed 

\vspace{4mm}

{\sc Examples.} 1. Let $X$ be an irreducible $k$-variety with the 
function field $L$ embedded into $F$ and $Q\subseteq{\mathbb Q}
[\{L\stackrel{/k}{\hookrightarrow}F\}]$ be a $G$-submodule such 
that ${\mathbb Q}[\{L\stackrel{/k}{\hookrightarrow}F\}]/Q$ is 
a quotient of $CH_0(X_F)_{{\mathbb Q}}$. Then it follows from the 
moving lemma of \cite{repr} and from Corollary \ref{cohom-dim} 
that $H^1_{{\mathcal S}m_G}(G_{F|F'},Q)=0$. 

2. (`Local acyclicity'). Suppose that the whole 
category ${\mathcal S}m_G$ satisfies the equivalent conditions of the 
previous lemma. Then for any finitely generated smooth representation 
$W$ of $G$ there exists an open subgroup $U\subset G$ such that 
$H^{>0}_{{\mathcal S}m_G}(U',W)=0$ for any open subgroup $U'\subseteq U$ 
of type $G_{F|L}$. (Indeed, a choice of $N$ generators of $W$ determines 
a surjection $\pi:{\mathbb Q}[\{K\stackrel{/k}{\hookrightarrow}F\}]^N
\longrightarrow W$ for a field extension of finite type $K|k$. Then 
for any field extension $L|K$ of finite type both ${\mathbb Q}
[\{K\stackrel{/k}{\hookrightarrow}F\}]^N$ and the kernel of $\pi$ are 
acyclic with respect to $H^{>0}_{{\mathcal S}m_G}(G_{F|L},-)$. \qed)

\vspace{4mm}

An application, we have in mind is, e.g., to the semi-simplicity 
of the modules of {\sl regular} differential forms and of similar 
representations. This is related to the coincidence of homological 
and numerical equivalence relations on algebraic cycles. Some 
known results are collected in Appendix \ref{differential-forms}.

\section{The dominant topology} \label{top}
Consider the category ${\mathfrak D}m_k$ of smooth $k$-morphisms 
of smooth $k$-schemes endowed with the pretopology, 
where the covers of an object $X$ are dominant morphisms to $X$.  

Given a presheaf ${\mathcal F}$ on ${\mathfrak D}m_k$ one can 
extend it to (the spectra of) the filtered unions ${\mathcal O}
=\lim\limits_{_A\longrightarrow\phantom{_A}}A$ of finitely generated 
smooth $k$-subalgebras $A$ by ${\mathcal F}({\mathcal O}):=
\lim\limits_{_A\longrightarrow\phantom{_A}}{\mathcal F}({\rm Spec}(A))$. 
Clearly, this is independent of presentation of ${\mathcal O}$ as a filtered 
union. Our main examples of ${\mathcal O}$ will be fields. 
\begin{lemma} \label{ekviv-st} The functor ${\mathfrak W}:{\mathcal F}
\mapsto{\mathcal F}(F)$ is an equivalence between the category of 
sheaves of sets (resp., of abelian groups etc) on ${\mathfrak D}m_k$ and 
the category of smooth $G$-sets (resp., of $G$-modules etc). \end{lemma}
{\it Proof.} Conversely, to define a quasi-inverse functor ${\mathfrak F}$ 
from the category of smooth $G$-sets to the category of sheaves of sets on 
${\mathfrak D}m_k$, fix an embedding over $k$ into $F$ of the function 
field of each connected component of each object of ${\mathfrak D}m_k$. 
Then to each $G$-set $W$ and an object $\coprod\limits_{\alpha}U_{\alpha}$, 
for some collection $(U_{\alpha})_{\alpha}$ of irreducible objects of 
${\mathfrak D}m_k$ we associate 
$\prod\limits_{\alpha}W^{G_{F|k(U_{\alpha})}}$. Note, that each 
homomorphism $K\stackrel{\tau}{\hookrightarrow}L$ of subfields in $F$ 
of finite type over $k$ extends to an element $\widetilde{\tau}\in G$, and 
the double coset $G_{F|L}\widetilde{\tau}G_{F|K}$ depends only on $\tau$, 
so the morphism $W^{G_{F|K}}\stackrel{\tau}{\longrightarrow}W^{G_{F|L}}$ 
is independent of a particular choice of $\widetilde{\tau}$. Then to 
a morphism $\coprod\limits_{\beta\in T}V_{\beta}
\stackrel{f}{\longrightarrow}\coprod\limits_{\alpha\in S}U_{\alpha}$, 
given by a collection $(k(U_{f(\beta)})
\stackrel{\tau_{\beta}}{\hookrightarrow}k(V_{\beta}))_{\beta}$ 
we associate the morphism \begin{equation}\label{uravnitel} 
\prod\limits_{\alpha\in S}W^{G_{F|k(U_{\alpha})}}
\stackrel{(\tau_{\beta})}{\longrightarrow}\prod\limits_{\alpha\in S}
\prod\limits_{\beta\in T_{\alpha}}W^{G_{F|k(V_{\beta})}}=
\prod\limits_{\beta\in T}W^{G_{F|k(V_{\beta})}}.\end{equation} 

The sheaf condition means that (i) the morphism 
(\ref{uravnitel}) is injective and (ii) it is the equalizer of 
\begin{equation}\label{dve-strelki} 
\prod\limits_{\alpha\in S}\prod\limits_{\beta\in T_{\alpha}}
W^{G_{F|k(V_{\beta})}}\rightrightarrows
\prod\limits_{\alpha\in S}\prod\limits_{\beta,\beta'\in T_{\alpha}}
\prod\limits_xW^{G_{F|k(x)}},\end{equation} 
where $x$ runs over irreducible components of 
$V_{\beta}\times_{U_{\alpha}}V_{\beta'}$, for any cover 
$\coprod\limits_{\beta\in T}V_{\beta}\stackrel{f}{\longrightarrow}
\coprod\limits_{\alpha\in S}U_{\alpha}$. The injectivity of 
(\ref{uravnitel}) is evident. To verify the condition (ii), 
we apply \cite[Proposition 2.14]{repr}: {\it if $L_1$ and $L_2$ are 
subextensions in $F|k$ of finite type such that the intersection of 
algebraic closures in $F$ of $L_1$ and $L_2$ is algebraic over 
$L_1\bigcap L_2$ then the subgroups $G_{F|L_1}$ and $G_{F|L_2}$ 
in $G$ generate $G_{F|L_1\bigcap L_2}$.} Namely, if 
$(w_{\beta})_{\beta\in T}$ belongs to the equalizer then 
$w_{\beta}\in W^{G_{F|k(U_{f(\beta)})}}$ and $w_{\beta}=w_{\beta'}$ 
if $f(\beta)=f(\beta')$, which means that $(w_{\beta})_{\beta\in T}$ 
belongs to the image of (\ref{uravnitel}).

Evidently, ${\mathfrak W}{\mathfrak F}\cong Id$. To check that 
${\mathfrak F}{\mathfrak W}\cong Id$, it suffices to verify that 
${\mathcal F}(U)=({\mathfrak W}{\mathcal F})^{G_{F|k(U)}}$. 
But this is exactly the sheaf condition for the covers 
of $U$ by single connected elements. \qed

\section{`Acyclicity' of certain restrictions of $G$-modules}
\begin{proposition} Suppose that $F'$ is an extension 
of $k$ in $F$ of infinite transcendence degree. 

Then $H^{>0}_{{\mathcal S}m_G}(G_{F|\overline{F'}},W)
=H^{>0}_{{\mathcal S}m_G}(G_{F|F'},W\otimes{\mathbb Q})=0$ 
for any smooth $G$-module $W$. \end{proposition}
{\it Proof.} Let $I^{\bullet}$ be an injective resolution of $W$ in the 
category ${\mathcal S}m_G$. As the functor $H^0(G_{F|\overline{F'}},-):
{\mathcal S}m_G\longrightarrow{\mathcal S}m_{G_{\overline{F'}|k}}$ is 
an equivalence of categories (cf. \cite{repr}), it is exact. This gives 
the vanishing of $H^{>0}_{{\mathcal S}m_G}(G_{F|\overline{F'}},W):=
H^{>0}\left((I^{\bullet})^{G_{F|\overline{F'}}}\right)$. 

As the functor $H^0(G_{F|F'},-):{\mathcal S}m_G({\mathbb Q})
\longrightarrow{\mathcal A}b$ is a composition of exact functors 
$H^0(G_{F|\overline{F'}},-):{\mathcal S}m_G({\mathbb Q})
\longrightarrow{\mathcal S}m_{G_{\overline{F'}|k}}({\mathbb Q})$ 
and $H^0(G_{\overline{F'}|F'},-):
{\mathcal S}m_{G_{\overline{F'}|k}}({\mathbb Q})\longrightarrow
{\mathcal A}b$, it is exact itself, and thus, 
$H^{>0}_{{\mathcal S}m_G}(G_{F|F'},W\otimes{\mathbb Q}):=
H^{>0}\left((I^{\bullet}\otimes{\mathbb Q})^{G_{F|F'}}\right)=0$. \qed 

\section{\v{C}ech cohomology}
\begin{lemma} \label{sta-abs} Let $S\stackrel{\pi}{\longrightarrow}T$ 
be a map of sets. For each $q\ge 0$, let $\stackrel{\circ}{S^q_T}\subseteq
\times^q_TS$ be a subset of the fibre product such that restriction to 
$\stackrel{\circ}{S^q_T}$ of each projection factors through 
$\stackrel{\circ}{S^{q-1}_T}$ and for any finite subset $U$ in 
$\stackrel{\circ}{S^q_T}$ over any element of $T$ (i.e. with 
$\#(\pi(U))=1$) there is $u\in S$ such that $\{u\}\times U\subseteq
\stackrel{\circ}{S^{q+1}_T}$. Then the natural complex (where $d_i$ 
are alternating sums of projections) 
${{\mathbb Q}[\stackrel{\circ}{S^0_T}]}\stackrel{d_1}{\longleftarrow}
{{\mathbb Q}[\stackrel{\circ}{S^1_T}]}\stackrel{d_2}{\longleftarrow}
{{\mathbb Q}[\stackrel{\circ}{S^2_T}]}\stackrel{d_3}{\longleftarrow}
{{\mathbb Q}[\stackrel{\circ}{S^3_T}]}\stackrel{d_4}{\longleftarrow}
\dots$ is acyclic. \end{lemma}
{\it Proof.} Let $\alpha=\sum_ia_i(s_{1i},\dots,s_{qi})\in\ker d_q$. 
Then \begin{multline*}d_{q+1}:\sum_ia_i(u_i,s_{1i},\dots,s_{qi})
\mapsto\sum_ia_i(s_{1i},\dots,s_{qi})+\sum_i\sum_{j=1}^q(-1)^j
a_i(u_i,s_{1i},\dots,\widehat{s_{ji}},\dots,s_{qi})=\alpha
\end{multline*} for any collection $\{u_i\}_i$ such that 
$(u_i,s_{1i},\dots,s_{qi})\in\stackrel{\circ}{S^{q+1}_T}$ for all 
$i$ and $u_i$ depends only on the projection of $s_{1i}$ to $T$. \qed 

\vspace{4mm}

\label{homotopy-invariant}
Denote by ${\mathcal I}_G$ the full subcategory of ${\mathcal S}m_G$, whose 
objects $W$ satisfy $W^{G_{F|L'}}=W^{G_{F|L}}$ for any purely transcendental 
extension $L'|L$ in $F|k$. More on this can be found in \cite{repr}. 
The equivalence of categories of Lemma \ref{ekviv-st} restricts to an 
equivalence of categories between ${\mathcal I}_G$  and the category 
of ${\mathbb A}^1$-{\sl invariant} sheaves on ${\mathfrak D}m_k$: 
${\mathcal F}(X)={\mathcal F}(X\times{\mathbb A}^1)$ for any smooth 
$k$-variety $X$. This is why the objects of ${\mathcal I}_G$ are called 
`homotopy invariant' representations. 

Let ${\mathcal I}:{\mathcal S}m_G\longrightarrow{\mathcal I}_G$ be the 
left adjoint of the inclusion functor ${\mathcal I}_G\hookrightarrow
{\mathcal S}m_G$ and $C_{k(Y)}:={\mathcal I}{\mathbb Q}[\{k(Y)
\stackrel{/k}{\hookrightarrow}F\}]$ for any irreducible $k$-variety $Y$. 
In the case of smooth proper $Y$, the object $C_{k(Y)}$ is related to 
the Chow group of zero-cycles on $Y$: there is a natural surjection 
$C_{k(Y)}\longrightarrow\hspace{-3mm}\to CH_0(Y_F)_{{\mathbb Q}}$. 
\begin{lemma} \label{nad-tochkoy} The complex 
$\check{C}^{L(Y)|L}_{\bullet}:=(\dots\to
{\mathbb Q}[\{L(Y^2)\stackrel{/k}{\hookrightarrow}F\}]\to
{\mathbb Q}[\{L(Y)\stackrel{/k}{\hookrightarrow}F\}]\to
{\mathbb Q}[\{L\stackrel{/k}{\hookrightarrow}F\}]\to 0)$ 
is acyclic for any field extension $L|k$ and any $L$-variety $Y$. 
If $Y$ is a smooth proper $k$-variety then the complexes 
$\dots\to C_{k(Y^2)}\to C_{k(Y)}\to{\mathbb Q}\to 0$ 
and $\dots\to CH_0(Y^2)_{{\mathbb Q}}\to CH_0(Y)_{{\mathbb Q}}\to
{\mathbb Q}\to 0$ are also acyclic. \end{lemma}
{\it Proof.} The acyclicity of $\check{C}^{L(Y)|L}_{\bullet}$ follows from 
Lemma \ref{sta-abs} applied to $T=\{L\stackrel{/k}{\hookrightarrow}F\}$ 
and $\stackrel{\circ}{S^q_T}=\{L(Y^q)\stackrel{/k}{\hookrightarrow}F\}$. 

Let $\alpha$ be an element of either $C_{k(Y^q)}$, or of 
$CH_0(Y^q)_{{\mathbb Q}}$ annihilated by the differential. As the 
morphisms ${\mathbb Q}[\{k(Y^q)\stackrel{/k}{\hookrightarrow}F\}]
\to C_{k(Y^q)}\to CH_0(Y^q_F)_{{\mathbb Q}}$ are surjective, $\alpha$ 
is represented by $\sum_ia_i(s_{1i},\dots,s_{qi})$. Let $u\in Y(F)$ 
be a generic point in generic position with respect to all $s_{ji}$, 
and $\beta$ be the image of $\sum_ia_i(u,s_{1i},\dots,s_{qi})$ 
in $C_{k(Y^q)}$, or in $CH_0(Y^q)_{{\mathbb Q}}$. Then 
$\partial\beta=\alpha+\sum_ia_i\sum_{j=1}^q(-1)^j
[(u,s_{1i},\dots,\widehat{s_{ji}},\dots,s_{qi})]$. 

It remains to show that $\sum_ia_i\sum_{j=1}^q(-1)^j
[(u,s_{1i},\dots,\widehat{s_{ji}},\dots,s_{qi})]=0$, 
which is clear in the case of $CH_0$. 

For the case of $C_{k(Y^{\bullet})}$ note that for any smooth 
representation $W$ the kernel of the surjection $W\to{\mathcal I}W$ 
is generated by the elements $\sigma w-w$ for $w\in W^{G_{F|L'}}$, 
$\sigma\in G_{F|L}$ and all purely transcendental extensions $L'|L$. 
Let $W={\mathbb Q}[\{k(Y^{q-1})\stackrel{/k}{\hookrightarrow}F\}]$. 
Then ${\mathcal I}W=C_{k(Y^{q-1})}$ and $\sum_ia_i
\sum_{j=1}^q(-1)^j(s_{1i},\dots,\widehat{s_{ji}},\dots,
s_{qi})=\sum_s(\sigma_sw_s-w_s)$, where $w_s\in W^{G_{F|L'_s}}$, 
$\sigma_s\in G_{F|L_s}$ and extensions $L'_s|L_s$ are purely 
transcendental and of finite type over $k$. 

Let $u:k(Y)\stackrel{/k}{\hookrightarrow}F$ be a generic point of $Y$ 
in generic position with respect to all $L'_s$ and $\sigma_s(L'_s)$. 
Then there exists a collection of $\sigma'_s\in G_{F|L_s}$ 
such that $\sigma'_sw_s=\sigma_sw_s$ and $\sigma'_su=u$, so 
$\sum_ia_i\sum_{j=1}^q(-1)^j(u,s_{1i},\dots,\widehat{s_{ji}},
\dots,s_{qi})=\sum_s(\sigma'_s(u\otimes w_s)-u\otimes w_s)\in
{\mathbb Q}[\{k(Y^q)\stackrel{/k}{\hookrightarrow}F\}]\subset 
{\mathbb Q}[\{k(Y)\stackrel{/k}{\hookrightarrow}F\}]\otimes W$, 
i.e. $\sum_ia_i\sum_{j=1}^q(-1)^j[(u,s_{1i},\dots,
\widehat{s_{ji}},\dots,s_{qi})]=0$ in $C_{k(Y^q)}$, since 
$u\otimes w_s$ is fixed by $G_{F|L'_su(k(Y))}$, $L'_su(k(Y))$ 
is purely transcendental over $L_su(k(Y))$ and $\sigma'_s\in 
G_{F|L_su(k(Y))}$. \qed

\begin{corollary} \label{tochnyy-cheh} For any $W\in{\mathcal I}_G$ 
the complex $0\to W^G\to W^{G_{F|k(Y)}}\to W^{G_{F|k(Y^2)}}\to\dots$ 
is exact. In particular, $0\to CH^q(X)_{{\mathbb Q}}\to 
CH^q(X_{k(Y)})_{{\mathbb Q}}\to CH^q(X_{k(Y^2)})_{{\mathbb Q}}
\to\dots$ is exact. \end{corollary}
{\it Proof.} The complex ${\rm Hom}_G(C_{k(Y^{\bullet})},W)$ 
is exact, since, by Lemma \ref{nad-tochkoy}, $C_{k(Y^{\bullet})}$ 
is a projective resolution of 0 in ${\mathcal I}_G$. \qed 

\begin{corollary}\label{tauto-sheaves} Sending the function field $L$ of 
a $k$-variety to ${\mathbb Q}[\{L\stackrel{/k}{\hookrightarrow}F\}]\in
{\mathcal S}m_G$, resp. to $C_L\in{\mathcal I}_G$, defines a sheaf on 
${\mathfrak D}m_k$ with values in ${\mathcal S}m_G^{{\rm op}}$, resp. in 
${\mathcal I}_G^{{\rm op}}$. {\rm (Since ${\mathcal I}$ is right exact.)} 
\qed \end{corollary}

Denote by $\check{H}^{\ast}(X,-)$ the \v{C}ech cohomology, i.e. 
$\lim\limits_{\longrightarrow}\check{H}^{\ast}_{{\mathfrak U}}(X,-)$ 
for ${\mathfrak U}$ running over the covers of $X$. 
\begin{corollary}[\cite{milne} Ch.III, Corollary 2.5] \label{crit}
$\check{H}^{\ast}$ coincides with $H^{\ast}$ for any sheaf if and only 
if $\check{H}^{\ast}$ transforms any short exact sequence of sheaves to 
a long exact sequence of \v{C}ech cohomologies. \end{corollary}

For any extension $K$ of $k$ in $F$ with $F$ of infinite 
transcendence degree over $K$ fix a transcendence basis 
$\{x_1,x_2,x_3,\dots\}=\{x^{(\overline{K})}_1,x^{(\overline{K})}_2,
x^{(\overline{K})}_3,\dots\}$ of $F$ over $K$. For each $m\ge 0$ 
set $F_m=\overline{K(x_{2^m},x_{2^m\cdot 3},x_{2^m\cdot 5},\dots)}$, 
and for each $j\ge 0$ fix a self-embedding $\sigma_j:F
\stackrel{/\overline{K}}{\hookrightarrow}F$ 
such that $\sigma_j|_{F_s}=id$ if $j>s$ and $\sigma_j|_{F_s}=
\sigma_s|_{F_s}:F_s\stackrel{\sim}{\longrightarrow}F_{s+1}$ if 
$j\le s$. For any extension $L$ of $K$ we fix $L_0\subseteq F_0$ 
isomorphic to $L$ over $K$ and set $L_s:=\sigma_0^s(L_0)$. 

For a $G$-module $W$ set $~_K^LW^{\bullet}:=(
W^{G_{F|L_0}}\stackrel{\sigma_0-\sigma_1}{-\!\!\!\longrightarrow}
W^{G_{F|L_0L_1}}\stackrel{\sigma_0-\sigma_1+\sigma_2}
{-\!\!\!-\!\!\!-\!\!\!-\!\!\!\longrightarrow}
W^{G_{F|L_0L_1L_2}}\longrightarrow\dots)$. 

\begin{proposition} \label{Chech=ext} One has 
$H^q_{{\mathcal S}m_G}(G_{F|K},W)=H^q(~_K^FW^{\bullet})$ for 
any smooth representation $W$ of $G$ and any algebraically 
closed extension $K\supseteq k$ in $F$. \end{proposition}
{\it Proof.} According to Corollary \ref{crit}, it suffices 
to check that the functor $W\longmapsto H^q(~_K^FW^{\bullet})$ 
coincides with the \v{C}ech cohomology $\lim\limits_{\longrightarrow}
H^q({\mathcal F}(U)\to{\mathcal F}(U\times_kU)\to{\mathcal F}
(U\times_kU\times_kU)\to\dots)$ for the sheaf ${\mathcal F}$ on 
${\mathfrak D}m_k$ corresponding to $W$, and transforms any short 
exact sequence of representations to a long exact sequence of spaces. 
By Lemma 5.7 of \cite{repr}, the functor $W\longmapsto~_K^FW^{\bullet}$ 
from ${\mathcal S}m_G$ to the category of complexes of 
${\mathbb Q}$-vector spaces is exact, so the latter property follows. 

If $K$ is algebraically closed in $L$ and $L$ is of finite type 
over $K$ then $_K^LW^{\bullet}$ is isomorphic to the \v{C}ech 
complex ${\rm Hom}_G(\check{C}^{L|K}_{\bullet},W)$ of the cover 
${\rm Spec}(L)\longrightarrow{\rm Spec}(K)$ with coefficients 
in ${\mathcal F}_W$. Therefore, under these assumptions, $W
\longmapsto H^q(~_K^LW^{\bullet})$ is the \v{C}ech cohomology 
of the cover ${\rm Spec}(L)\longrightarrow{\rm Spec}(K)$. 
As $_{\overline{K}}^FW^{\bullet}=
\lim\limits_{_L\longrightarrow\phantom{_L}}~_{\overline{K}}^LW^{\bullet}$, 
we see that $H^q(~_{\overline{K}}^FW^{\bullet})=\lim\limits
_{_L\longrightarrow\phantom{_L}}H^q(~_{\overline{K}}^LW^{\bullet})$, 
where $L$ runs over the set of subfields in $F$ of finite type 
over $\overline{K}$. And finally, the functor $W\longmapsto 
H^q(~_{\overline{K}}^FW^{\bullet})$ coincides with the 
\v{C}ech cohomology of ${\rm Spec}(\overline{K})$. \qed

\vspace{4mm}

{\sc Examples.} 1. Let $W$ be a smooth $G$-module with trivial 
restriction to $G_{F|\overline{k}}$. Then ${\mathcal F}_W=\pi^{-1}
{\mathcal F}_W|_{{\rm Spec}(k)_{{\rm \acute{e}t}}}$, where 
${\mathfrak D}m_k\stackrel{\pi}{\longrightarrow}
{\rm Spec}(k)_{{\rm \acute{e}t}}$ sends a smooth $k$-variety $U$ 
to the maximal $k$-scheme $U'$ \'{e}tale over ${\rm Spec}(k)$ in the 
Stein decomposition $U\longrightarrow U'\longrightarrow{\rm Spec}(k)$. 

For any connected cover $U\longrightarrow{\rm Spec}(\overline{k})$ 
the \v{C}ech complex with coefficients in ${\mathcal F}_W$ is just 
$W\stackrel{0}{\longrightarrow}W\stackrel{id}{\longrightarrow}W
\stackrel{0}{\longrightarrow}W\stackrel{id}{\longrightarrow}\dots$, so 
$\check{H}^{>0}_{{\mathfrak D}m_{\overline{k}}}({\mathcal F}_W)=0$. 

Then the Hochschild--Serre spectral sequence $E^{p,q}_2=
H^p(G_{\overline{k}|k},H^q(G_{F|\overline{k}},W))$ is degenerate: \\
$E^{p,>0}_2=0$, and thus, $\check{H}^{\ast}_{{\mathfrak D}m_k}
({\mathcal F}_W)=E^{\ast,0}_2=H^{\ast}_{{\rm \acute{e}t}}
({\rm Spec}(k),{\mathcal F}_W|_{{\rm Spec}(k)_{{\rm \acute{e}t}}})$. 

2. The sheaf $L\mapsto{\mathbb Q}[\{L
\stackrel{/k}{\hookrightarrow}F\}]$ (see Corollary \ref{tauto-sheaves}) 
with values in ${\mathcal S}m_G^{{\rm op}}$ is acyclic.

3. Let $A$ be an abelian variety over $k$. Then for any smooth variety 
$U$ over $k$ and its smooth compactification $\overline{U}$ the 
\v{C}ech complex $B\stackrel{d_1}{\longrightarrow}B^2\stackrel{d_2}
{\longrightarrow}B^3\stackrel{d_3}{\longrightarrow}\dots$ of 
$A(F)/A(k)$, where $B^q={\rm Hom}({\rm Alb}(\overline{U}^q),A)$, is 
acyclic, since $$d_q:(b_1,\dots,b_q)\mapsto\left\{\begin{array}{ll} 
(0,b_1+b_2,0,b_3+b_4,0,\dots,0,b_{q-1}+b_q,0) & \mbox{if $q$ is even}\\ 
(-b_1,b_1,-b_3,b_3,\dots,-b_q,b_q) & \mbox{if $q$ is odd} 
\end{array}\right.$$

More generally, fix $n\ge 0$ and set $C^q=H^n(\overline{U}^q)/N^1$. 
The complex $C^{\bullet}$ is the direct sum of the sub-complexes 
$C^{\bullet}_{m_1,\dots,m_s}$ for all collections $(m_1,\dots,m_s)$ 
such that $m_1,\dots,m_s\ge 1$ and $m_1+\dots+m_s=n$, where 
$C^q_{m_1,\dots,m_s}=\!\!\!\!\!\!\!\bigoplus\limits_{i_0,\dots,i_s\ge 
0,~~i_0+\dots+i_s=q-s}\!\!\!\!\!\!\!H^0(\overline{U})^{\otimes i_0}
\otimes H^{m_1}(\overline{U})\otimes H^0(\overline{U})^{\otimes i_1}
\otimes H^{m_2}(\overline{U})\otimes\dots\otimes 
H^{m_s}(\overline{U}^q)\otimes H^0(\overline{U})^{\otimes i_s}/N^1
\cong(C^s_{m_1,\dots,m_s})^{\binom{q}{s}}$. 
Clearly, all $C^{\bullet}_{m_1,\dots,m_s}$ are acyclic. 

\subsection{`Geometric' representations and the acyclicity} 
Any field extension $K$ of $k$ is a filtered union of finitely 
generated smooth $k$-subalgebras. Then the construction of \S\ref{top}, 
preceding Lemma \ref{ekviv-st}, associates to any presheaf ${\mathcal F}$ 
on ${\mathfrak D}m_k$ the smooth ${\rm Aut}(K|k)$-set ${\mathcal F}(K)$.
\begin{lemma} One has ${\mathcal F}(F)^{G_{F|F'}}={\mathcal F}(F')$ 
for any $F'=\overline{F'}\subseteq F$ with ${\rm tr.deg}(F'|k)=\infty$. 
\end{lemma}
(As it shows the example of the presheaf ${\mathcal F}:U\mapsto
\Gamma(\overline{U},\bigotimes^2_{{\mathcal O}}\Omega^1_{\overline{U}/k})$, 
one cannot drop the condition that $F'$ is algebraically closed. 
Here $\overline{U}$ is a smooth compactification of $U$.) 

{\it Proof.} Fix an isomorphism $\alpha:F\stackrel{\sim}{\longrightarrow}F'$ 
over $k$. Then $\alpha^{\ast}:\lim\limits_{_U\longrightarrow\phantom{_U}}
{\mathcal F}(U)\stackrel{\sim}{\longrightarrow}\lim\limits
_{_V\longrightarrow\phantom{_V}}{\mathcal F}(V)$, 
where ${\mathcal O}(U)\subset F$ and ${\mathcal O}(V)\subset F'$ 
are smooth. For any $U$ there is $\sigma\in G$ such that $\sigma|
_{{\mathcal O}(U)}=\alpha|_{{\mathcal O}(U)}$, so $\lim\limits
_{_V\longrightarrow\phantom{_V}}{\mathcal F}(V)=\lim
\limits_{_{(U,\sigma)}\longrightarrow\phantom{_{(U,\sigma)}}}\sigma
{\mathcal F}(U)={\mathcal F}(F)^{G_{F|F'}}$. \qed 

\vspace{4mm}

Note that the category of representations of type ${\mathcal F}(F)$ 
for ${\mathcal F}$ taking values in commutative groups is tensor. 
More generally, given a presheaf of commutative rings ${\mathcal A}$, 
the representations ${\mathcal F}(F)$ for the ${\mathcal A}$-modules 
${\mathcal F}$ are ${\mathcal A}(F)$-modules and they form a tensor 
category with respect to the operation $\otimes_{{\mathcal A}(F)}$. 

Assume that a presheaf ${\mathcal F}$ on ${\mathfrak D}m_k$ is endowed 
with transformations $i^x_{X,Y}:{\mathcal F}(X\times_kY)\longrightarrow
{\mathcal F}(Y)$ for any smooth $X,Y$ and any $x\in X(k)$ such that 
$i^x_{X,Y}\circ{}_X{\rm pr}^{\ast}_Y=id_{{\mathcal F}(Y)}$ and 
$i^x_{X,Y\times_kZ}\circ{}_Z{\rm pr}^{\ast}_{X\times_kY}=
{}_Z{\rm pr}^{\ast}_Y\circ i^x_{X,Y}$, where 
$_X{\rm pr}_Y:X\times_kY\longrightarrow Y$ is the projection. 
\begin{corollary}\label{2-3-acikl} $H^{>0}(G,W)=0$ for any 
$W\in{\mathcal I}_G$; $H^1(G,{\mathcal F}(F))=\{\ast\}$ for any 
group-valued ${\mathcal F}$; and $H^{>0}(G,{\mathcal F}(F))=0$ for 
any ${\mathcal F}$ with values in commutative groups. \end{corollary}
{\it Proof.} By Proposition \ref{Chech=ext}, $\check{H}^{\ast}(G,-)
=H^{\ast}(G,-)$. By Corollary \ref{tochnyy-cheh}, 
$\check{H}^{\ast}(G,W)=W^G$. 

It suffices to show that the complex $0\to{\mathcal F}(k)\stackrel
{d_0}{\longrightarrow}{\mathcal F}(X)\stackrel{d_1}{\longrightarrow}
{\mathcal F}(X^2)\stackrel{d_2}{\longrightarrow}{\mathcal F}(X^3)
\stackrel{d_3}{\longrightarrow}\dots$ is exact for any smooth 
$k$-variety $X$, where $d_q=\sum_{i=0}^q(-1)^i({\rm pr}_q^i)^{\ast}$ 
and ${\rm pr}_q^i:X^{q+1}\longrightarrow X^q$ omits the $i$-th
multiple. Fix some $u\in X(k)$ and denote by $u_q$ the closed 
embedding $X^q\stackrel{u\times}{\hookrightarrow}X^{q+1}$. 
Then ${\rm pr}_q^i\circ u_q=u_{q-1}\circ{\rm pr}_{q-1}^{i-1}$ 
for any $1\le i\le q$ and ${\rm pr}_q^0\circ u_q=id_{X^q}$, 
so $u_q^{\ast}\circ d_q+d_{q-1}\circ u_{q-1}^{\ast}=
id_{{\mathcal F}(X^q)}$, where $u_q^{\ast}:=i^u_{X,X^q}$, i.e. 
$u_{\bullet}^{\ast}$ defines a contracting homotopy. \qed 

\vspace{4mm}

{\sc Examples of} ${\mathcal F}(F)$ are $\bigotimes^{\bullet}_F
\Omega^1_{F|k_0}$, $\bigotimes^{\bullet}_{k_0}\Omega^1_{F|k_0,{\rm cl}}$, 
$\bigotimes^{\bullet}_{k_0}\Omega^1_{F|k_0,{\rm ex}}$ for any 
$k_0\subseteq k$ and $A(F)$ for any group $k$-variety $A$, or 
$A(F)/A(k)$ if $A$ is commutative. More examples are in the following 

\begin{corollary} \label{cohom-dim} Let $E$ be a field, endowed with 
a smooth $G$-action, which is presentable as ${\mathcal G}(F)$ for 
a presheaf of commutative rings ${\mathcal G}$ on ${\mathfrak D}m_k$ 
admitting transformations $i^x_{X,Y}$ as above. 
\begin{enumerate} 
\item \label{acyclich-obrazuyuxie} Let $B$ be such a family of fields 
of finite type over $k$ that any extension of $k$ of finite type can 
be embedded (over $k$) into an element of $B$. 

Then the collection 
$\{E[\{L\stackrel{/k}{\hookrightarrow}F\}]\}_{L\in B}$ is a system of 
acyclic generators of the category ${\mathcal S}m_G(E)$ (of 
smooth $E$-semilinear representations of $G$). 
\item The cohomological dimension of the category 
${\mathcal S}m_G(E)$ is infinite. \end{enumerate} \end{corollary}
{\it Proof.} \begin{enumerate} \item For any field extension 
$L|k$ of finite type consider the functor ${\mathcal F}(X)=Z^q(X_L)
\otimes{\mathcal G}(X)$, where $q={\rm tr.deg}(L|k)$. 
Then $E[\{L\stackrel{/k}{\hookrightarrow}F\}]={\mathcal F}(F)$. 
Define the transformations $i^x_{X,Y}:Z^q(X\times_kY_L)\otimes
{\mathcal G}(X\times_kY)\longrightarrow Z^q(Y_L)\otimes{\mathcal G}(Y)$ 
by $\alpha\otimes g\mapsto\alpha\cap(\{x\}\times Y)\otimes i^x_{X,Y}g$ 
if the intersection is transversal, and $\alpha\otimes g\mapsto 0$ 
otherwise, for any irreducible $\alpha\subset X\times_kY_L$. 
\item For any $q\ge 1$ and any object $W\in{\mathcal S}m_G(E)$ such that 
$H^q(G,W)\neq 0$ there is a surjection $\alpha:U\longrightarrow W$, where 
$U$ is a direct sum of acyclic generators of ${\mathcal S}m_G(E)$, cf. 
(\ref{acyclich-obrazuyuxie}). Then there is an embedding 
$H^q(G,W)\hookrightarrow H^{q+1}(G,\ker\alpha)$, i.e., 
$H^{q+1}(G,W')\neq 0$ for $W'=\ker\alpha$. \qed \end{enumerate} 

\vspace{4mm}

{\sc Examples.} 1. $E=F={\mathcal G}_1(F)$ with 
${\mathcal G}_1(X)={\mathcal O}(X)$ and 
$i^x_{X,Y}$ being the restriction to $\{x\}\times Y$. 

2. ${\mathbb Q}\subseteq E={\mathcal G}_2(F)\subseteq k$ with 
${\mathcal G}_2:X\mapsto E^{\pi_0(X)}$ being the sub-presheaf of $E$-valued 
locally constant functions in ${\mathcal G}_1$ and $i^x_{X,Y}$ being the 
restriction of the corresponding transformation on ${\mathcal G}_1$. 

\begin{conjecture} ${\rm Ext}^{\ast}_{{\mathcal S}m_G}
(W_1,W_2)={\rm Ext}^{\ast}_{{\mathcal I}_G}(W_1,W_2)$ 
for any $W_1,W_2\in{\mathcal I}_G$. \end{conjecture}
This is known for $\ast\le 1$ (since ${\mathcal I}_G$ is a Serre subcategory 
of ${\mathcal S}m_G$) and for trivial $W_1$ (by Corollary \ref{2-3-acikl}). 

\subsection{${\mathbb A}^1$-invariance of some presheaves} 
\label{A-1-invar-puchki}
Let ${\mathcal V}_k$ be a category of $k$-varieties, 
containing all smooth varieties. Let ${\mathcal L}$ be a category, 
where all self-embeddings are isomorphisms (e.g., an abelian category 
such that for any object the multiplicities of its irreducible subquotients 
are finite\footnote{The multiplicity of an irreducible object $X$ in $W$ is 
defined inductively: it is 0 if for any filtration $W\supseteq Y\supseteq Z$ 
the quotient $Y/Z$ is not isomorphic to $X$; it is $N>0$ if there is 
a filtration $W\supseteq Y\supseteq Z$ such that $Y/Z\cong X$ and 
the sum the multiplicities of $X$ in $W/Y$ and in $Z$ is $N-1$. By 
Jordan--H\"{o}lder theorem, the multiplicity is well-defined.}). 

An ${\mathcal L}$-valued presheaf ${\mathcal F}$ on ${\mathcal V}_k$ 
is ${\mathbb A}^1$-{\sl invariant}, if ${\mathcal F}(U)=
{\mathcal F}(U\times{\mathbb A}^1)$ for any $U\in{\mathcal V}_k$. 

Consider any pretopology on ${\mathcal V}_k$ such that 
${\mathbb A}^1_k\longrightarrow{\rm Spec}(k)$ is a covering 
and one of the following two conditions holds:
\begin{itemize} \item \label{opred-dopust} $k$ contains all $2$-primary roots 
of unity and there are no non-trivial divisible involutions in the group 
${\rm Aut}(C)$ for any $C\in{\mathcal L}$ (e.g., if ${\mathcal L}$ is the 
category of finitely generated admissible\footnote{A smooth representation of 
a totally disconnected topological group is called {\sl admissible} if any 
open subgroup fixes a finite-dimensional subspace in it.} representations 
over a number field of a totally disconnected topological group); 
\item $({\mathbb A}^1_k\smallsetminus\{0\})\coprod({\mathbb A}^1_k
\smallsetminus\{1\})\stackrel{id\sqcup id}{\longrightarrow}{\mathbb A}^1_k$ 
and ${\mathbb G}_{m,k}\longrightarrow{\mathbb G}_{m,k}$, $x\mapsto x^2$, 
are coverings (in particular, ${\mathcal F}({\mathbb A}^1_X)\longrightarrow
{\mathcal F}({\mathbb G}_{m,X})$ is injective for any $X\in{\mathcal V}_k$ 
and any sheaf on ${\mathcal V}_k$). 
\end{itemize} 

\begin{proposition} \label{A1-invar} Let ${\mathcal F}$ be a sheaf on 
${\mathcal V}_k$ with values in ${\mathcal L}$. Then ${\mathcal F}$ 
is ${\mathbb A}^1$-invariant. \end{proposition}
{\it Proof.} Let $U\in{\mathcal V}_k$ and $\sigma$ be the involution of 
${\mathcal F}(U\times{\mathbb A}^1\times{\mathbb A}^1)$ induced by the 
interchanging of the two multiples 
${\mathbb A}^1$, so ${\rm pr}_1^{\ast}=\sigma\circ{\rm pr}^{\ast}_2$, 
where ${\rm pr_1},{\rm pr}_2:~U\times{\mathbb A}^1\times{\mathbb A}^1
\rightrightarrows U\times{\mathbb A}^1$ are the projections. 

As $U\times{\mathbb A}^1\longrightarrow U$ is a covering, we get 
that ${\mathcal F}(U)$ is the equalizer of the injections 
$${\rm pr}_1^{\ast},{\rm pr}^{\ast}_2:~{\mathcal F}(U\times{\mathbb A}^1)
\rightrightarrows{\mathcal F}(U\times{\mathbb A}^1\times{\mathbb A}^1).$$ 

If $k$ contains all $2$-primary roots of unity then $\sigma$ 
is a divisible involution in the group ${\rm Aut}_{{\mathcal L}}
({\mathcal F}(U\times{\mathbb A}^1\times{\mathbb A}^1))$. 

If the Zariski covering ${\mathbb G}_{m,k}\coprod{\mathbb G}_{m,k}
\longrightarrow{\mathbb A}^1_k$ is a covering in ${\mathcal V}_k$ then 
${\mathcal F}(U\times{\mathbb A}^1\times{\mathbb A}^1)\longrightarrow
{\mathcal F}(U\times(({\mathbb A}^1\times{\mathbb A}^1)\smallsetminus
\Delta_{{\mathbb A}^1}))$ is injective, since $({\mathbb A}^1\times
{\mathbb A}^1)\smallsetminus\Delta_{{\mathbb A}^1}\stackrel{\sim}
{\longrightarrow}{\mathbb A}^1\times{\mathbb G}_m$, $(x,y)\mapsto(x,x-y)$. 

The morphism $U\times(({\mathbb A}^1\times{\mathbb A}^1)\smallsetminus
\Delta_{{\mathbb A}^1})\longrightarrow U\times(({\mathbb A}^1\times
{\mathbb A}^1)\smallsetminus\Delta_{{\mathbb A}^1})/{\mathfrak S}_2$ 
in the category ${\mathcal V}_k$ induces a ${\mathfrak S}_2$-invariant 
morphism ${\mathcal F}(U\times(({\mathbb A}^1\times{\mathbb A}^1)
\smallsetminus\Delta_{{\mathbb A}^1})/{\mathfrak S}_2)\longrightarrow
{\mathcal F}(U\times(({\mathbb A}^1\times{\mathbb A}^1)\smallsetminus
\Delta_{{\mathbb A}^1}))$ in ${\mathcal L}$, which is injective if 
${\mathbb G}_{m,k}\longrightarrow{\mathbb G}_{m,k}$, $x\mapsto x^2$, 
is a covering. As $(({\mathbb A}^1\times{\mathbb A}^1)\smallsetminus
\Delta_{{\mathbb A}^1})/{\mathfrak S}_2\stackrel{\sim}{\longrightarrow}
({\mathbb A}^1\times{\mathbb A}^1)\smallsetminus\Delta_{{\mathbb A}^1}$, 
$(x,y)\mapsto(x+y,x+y+(x-y)^2)$, the involution $\sigma$ acts trivially 
on ${\mathcal F}(U\times(({\mathbb A}^1\times{\mathbb A}^1)\smallsetminus
\Delta_{{\mathbb A}^1}))$ by the assumptions on ${\mathcal L}$. 

In any case, one has ${\rm pr}_1^{\ast}=\sigma\circ{\rm pr}^{\ast}_2=
{\rm pr}^{\ast}_2$, and thus, ${\mathcal F}(U)\longrightarrow
{\mathcal F}(U\times{\mathbb A}^1)$ is an isomorphism. \qed 

\vspace{4mm}

{\it Example.} If a morphism $p:U\longrightarrow Y$ in ${\mathcal V}_k$ 
admits a section $s:Y\longrightarrow U$ (e.g. for $U=Y\times_kZ$, and $p$ 
being the projection to $Y$, where $Z$ has a rational $k$-point) then the 
sequence ${\mathcal G}(Y)\stackrel{p^{\ast}}{\longrightarrow}{\mathcal G}(U)
\rightrightarrows{\mathcal G}(U\times_YU)$ is exact for any presheaf 
${\mathcal G}$: $p^{\ast}$ is injective, since the identity composition 
$Y\longrightarrow U\longrightarrow Y$ induces identity composition 
${\mathcal G}(Y)\longrightarrow{\mathcal G}(U)\longrightarrow{\mathcal G}(Y)$; 
$p^{\ast}$ is an equalizer of ${\mathcal G}(U)\rightrightarrows{\mathcal G}
(U\times_YU)$, since the compositions $U\stackrel{(id,sp)}{\longrightarrow}
U\times_YU\rightrightarrows U$ induce the identical and $(sp)^{\ast}=
p^{\ast}s^{\ast}$ compositions ${\mathcal G}(U)\rightrightarrows
{\mathcal G}(U\times_YU)\longrightarrow{\mathcal G}(U)$, i.e. any 
element in the equalizer belongs to the image of $p^{\ast}$. 

{\it Question.} Any cohomology theory $H^{\ast}$ on ${\mathcal V}_k$ with 
coefficients in a field and any $q\ge 0$ define an ${\mathbb A}^1$-invariant 
presheaf $H^q$ on ${\mathcal V}_k$. We want it to be a sheaf in some topology. 
Can one simply define the covering morphisms as smooth morphisms 
$p:U\longrightarrow X$ such that the sequence $0\longrightarrow 
H^{\ast}(Y)\longrightarrow H^{\ast}(U')\longrightarrow 
H^{\ast}(U'\times_YU')$ is exact for any base change 
$U'\longrightarrow Y$ of $p$?\footnote{It is not clear, whether this 
is compatible with compositions.} In particular, no non-surjective 
morphism in ${\mathcal V}_k$ is covering, if one allows the closed 
embeddings as morphisms in ${\mathcal V}_k$. (If $Y$ is a point of $X$, 
which is not covered, then $U'$ is empty, so the sequence $0\longrightarrow 
H^0(Y)\longrightarrow H^0(U')\longrightarrow H^0(U'\times_YU')$ is not exact.) 
It follows that this is the case for arbitrary morphisms admitting a section 
$U\longrightarrow X$.\footnote{and clearly this type of coverings is stable 
under base changes and compositions} Clearly, the Zariski coverings are not 
coverings for such topology (${\mathbb P}^n$ is a union of $n+1$ copies of 
${\mathbb A}^n$, but 
$H^{>0}({\mathbb P}^n)\longrightarrow H^{>0}({\mathbb A}^n)^{n+1}=0$ 
is not injective); the fibration ${\mathbb A}^{n+1}\smallsetminus\{0\}
\longrightarrow{\mathbb P}^n$ is not covering, since $H^{2n}({\mathbb P}^n)
\longrightarrow H^{2n}({\mathbb A}^{n+1}\smallsetminus\{0\})=0$ 
is not injective. 

\section{The site $\mathfrak{H}$ (the Hartogs topology)} \label{geom-puchki} 
Let $v:F^{\times}/k^{\times}\longrightarrow\hspace{-3mm}\to{\mathbb Q}^r$ 
be a discrete valuation, i.e. $F$ is algebraic over the subfield generated 
by an isomorphic lifting of the residue field of $v$ and by an isomorphic 
lifting of the valuation group ${\mathbb Q}^r$. Let $G_v$ be the decomposition 
subgroup of $v$ in $G$, i.e., the stabilizer of the valuation $k$-algebra 
${\mathcal O}_v\subset F$. 
For any $W\in{\mathcal S}m_G$ set $W_v:=\sum_{\sigma\in G_v}
W^{G_{F|\sigma(F')}}\subseteq W$. The additive functor 
${\mathcal S}m_G\longrightarrow{\mathcal S}m_{G_v}$, $W\mapsto W_v$, 
is fully faithful and preserves surjections and injections, cf. \cite[\S4.1]{max}. 

Let $\Gamma_r(W)$ be the intersection of $W_v$ over all discrete valuations 
$v:F^{\times}/k^{\times}\longrightarrow\hspace{-3mm}\to{\mathbb Q}^r$. Set 
$\Gamma:=\Gamma_1$, so $\Gamma_r:{\mathcal S}m_G\longrightarrow
{\mathcal S}m_G$ are additive functors. 

{\sc Examples.} ${\mathbb Q}[\{L\stackrel{/k}{\hookrightarrow}F\}]_v
={\mathbb Q}[\{L\stackrel{/k}{\hookrightarrow}{\mathcal O}_v\}]$; 
$F[\{L\stackrel{/k}{\hookrightarrow}F\}]_v
={\mathcal O}_v[\{L\stackrel{/k}{\hookrightarrow}{\mathcal O}_v\}]$. 

\begin{proposition}[\cite{max}] \label{loc-i-glob-invar} 
$\Gamma(W)=W_v=W$ for any $W\in{\mathcal I}_G$. \end{proposition}

\begin{lemma}[\cite{max}] \label{basic-glob} One has $(W_1\otimes W_2)_v
\subseteq(W_1)_v\otimes(W_2)_v$ and
$\Gamma(W_1\otimes W_2)\subseteq\Gamma(W_1)\otimes\Gamma(W_2)$ 
for any $W_1,W_2\in{\mathcal S}m_G$. However, $(W\otimes W)_v\neq 
W_v\otimes W_v$ if $W={\mathbb Q}[F\smallsetminus k]$. 

If either $W_1$ is a quotient of $A(F)$ for a commutative algebraic 
$k$-group $A$, or $W_1\in{\mathcal I}_G$ then $(W_1\otimes W_2)_v=
(W_1)_v\otimes(W_2)_v$ for any $W_2\in{\mathcal S}m_G$. \end{lemma} 

{\sc Remarks.} 1. By Proposition \ref{loc-i-glob-invar} and Lemma 
\ref{basic-glob}, $\Gamma(W\otimes F)=W\otimes k$ for any $W\in
{\mathcal I}_G$, so $\Gamma(V)\otimes F\longrightarrow\hspace{-3mm}\to V$ 
for any semilinear quotient $V$ of $W\otimes F$. 

2. If $W$ is an $F$-vector space then the 
$F$-vector space structure $F\otimes W\longrightarrow W$ on $W$ 
induces an ${\mathcal O}_v$-module structure $(F\otimes W)_v=
{\mathcal O}_v\otimes W_v\longrightarrow W_v$ on $W_v$. 
Clearly, $F\otimes_{{\mathcal O}_v}W_v\longrightarrow W$ is 
injective, but not surjective, as shows the example of 
$W=F[\{L\stackrel{/k}{\hookrightarrow}F\}]$. 

3. Clearly, $\Gamma_r$ preserves the injections, but not the 
surjections. Namely, let $W:=\bigotimes^N_kF\longrightarrow
\Omega^{N-1}_{F|k}$ be given by $a_1\otimes\dots\otimes a_N
\mapsto a_1da_2\wedge\dots\wedge da_N$. By Lemma \ref{basic-glob}, 
$W_v=\bigotimes^N_k{\mathcal O}_v$, so $(\bigotimes^{N-1}_F
\Omega^1_{F|k_0})_v=\bigotimes^{N-1}_F\Omega^1_{{\mathcal O}_v|k_0}$ 
for any $k_0\subseteq k$; and $\Gamma(\bigotimes_k^NF)=k$, but 
$\Gamma_r(\Omega^{\bullet}_{F|k})=\Omega^{\bullet}_{F|k,{\rm reg}}$ 
for any $r\ge 1$. 

\vspace{4mm}

For an integral normal $k$-variety $X$ with $k(X)\subset F$ 
let ${\mathfrak V}(X)$ be the set of all discrete valuations 
of $F$ of rank one trivial on $k$ such that their restrictions 
to $k(X)$ are either trivial, or correspond to divisors on $X$. 
Set ${\mathcal W}(X):=W^{G_{F|k(X)}}\cap
\bigcap_{v\in{\mathfrak V}(X)}W_v\subseteq W$. 

{\sc Remark.} $W^{G_{F|k(X)}}\cap W_v$ depends only on the restriction 
of $v$ to $k(X)$, since the set of $G_{F|k(X)}$-orbits $G_{F|k(X)}
\backslash G/G_v$ of the valuations of $F$ coincides with the set 
of valuations of $k(X)$ of rank $\le r$. E.g., if the restriction 
of $v$ to $k(X)$ is trivial then $W^{G_{F|k(X)}}\subseteq W_v$. 

Consider the following site $\mathfrak{H}$. Objects of 
$\mathfrak{H}$ are the smooth varieties over $k$. Morphisms in 
$\mathfrak{H}$ are the locally dominant morphisms, transforming 
non-dominant divisors to divisors. Coverings are smooth morphisms 
surjective over the generic point of any divisor on the target. 

\begin{lemma} A choice of $k$-embeddings into $F$ of all generic 
points of all smooth $k$-varieties defines a sheaf ${\mathcal W}$ 
on $\mathfrak{H}$ for any $W\in{\mathcal S}m_G$. \end{lemma} 
{\it Proof.} Clearly, if a dominant morphism $f:U\longrightarrow X$ 
transforms divisors on $U$, non-dominant over $X$, to divisors 
on $X$ then ${\mathfrak V}(U)\subseteq{\mathfrak V}(X)$, 
so ${\mathcal W}(X)\subseteq{\mathcal W}(U)$. 

If, moreover, the pull-back of any divisor on $X$ is 
a divisor on $U$ then ${\mathfrak V}(X)={\mathfrak V}(U)$. 

By Lemma \ref{ekviv-st}, $X\mapsto\prod_{x\in X^0}W^{G_{F|k(x)}}$ 
is a sheaf on ${\mathfrak D}m_k$, so the sequence $$0\longrightarrow
\prod_{x\in X^0}W^{G_{F|k(x)}}\longrightarrow\prod_{x\in U^0}
W^{G_{F|k(x)}}\longrightarrow\prod_{x\in(U\times_XU)^0}
W^{G_{F|k(x)}}$$ is exact. As ${\mathfrak V}(X)={\mathfrak V}(U)
={\mathfrak V}(U\times_XU)$, the sheaf property for the covering 
$f$ amounts to the exactness of the above sequence restricted 
to $\prod_{x\in U^0}\bigcap_{v\in{\mathfrak V}(X)}W_v$. \qed 

{\sc Examples.} $X\mapsto\Gamma(X,\Omega^q_{X|k})$ and 
$X\mapsto\Gamma(\overline{X},\Omega^q_{\overline{X}|k})$ are sheaves on 
$\mathfrak{H}$ for any integer $q\ge 0$; $X\mapsto\Gamma(\overline{X},
{\rm Sym}^2\Omega^q_{\overline{X}|k})$ is not a sheaf on $\mathfrak{H}$. 
Here $\overline{X}$ is a smooth compactification of $X$. 

\appendix
\section{Totally disconnected groups and their representations}
\label{general-topology}
Let $H$ be an arbitrary totally disconnected topological group, and $B$ 
a base of its open subgroups. By definition, this means that $B$ is a 
collection of subgroups of $H$ such that (i) it has the trivial 
intersection, (ii) the conjugate of each subgroup in $B$ by any element 
of $H$ contains a subgroup in $B$, (iii) each finite intersections of 
subgroups in $B$ contains a subgroup in $B$. Then a subgroup of $H$ is 
called {\sl open} if it contains an element of $B$ (so $H$ is Hausdorff). 

Let $E$ be a field of 
characteristic zero endowed with a smooth (e.g., trivial) $H$-action. 

An example of $H$ is given by a permutation group of a set (with the 
stabilizers of finite subsets as a base of open subgroups). This example 
is typical in the sense that arbitrary $H$ is a permutation group of the 
disjoint union $\coprod_{U\in B}H/U$ for arbitrary collection $B$ of 
subgroups of $H$ with trivial intersection of their conjugates, e.g., 
for a base of open subgroups of any totally disconnected group $H$. If 
there is an open subgroup $U$ of $H$ containing no nontrivial normal 
subgroups of $H$ then $H$ is a permutation group of the set $H/U$. 

In this appendix we make some general remarks on the category 
${\mathcal S}m_H(E)$ of $E$-vector spaces endowed with 
a smooth semilinear $H$-action, generalizing some well-known 
facts about locally compact totally disconnected groups. 

\subsection{Smooth $H$-sets and sheaves} \label{general-constr-sheaves}
It is well-known (e.g., \cite[Expos\'{e} IV, \S2.4--2.5]{sga-4-1} 
or \cite[\S8.1, Example 8.15 (iii)]{topos}) that the 
smooth $H$-sets and their $H$-equivariant maps form a topos.  

Let ${\mathfrak T}={\mathfrak T}(H,B)$ be the category whose objects 
are the elements of some base $B$ of open subgroups of $H$ and 
${\rm Hom}_{{\mathfrak T}}(U,V)=\{h\in H~|~hUh^{-1}\supseteq V\}/U$. 
The composition is defined in the natural way: if 
$U\stackrel{\alpha}{\longrightarrow}V\stackrel{\beta}{\longrightarrow}W$, 
i.e., $\alpha U\alpha^{-1}\supseteq V$ and $\beta V\beta^{-1}\supseteq W$, 
then $\beta\alpha U(\beta\alpha)^{-1}\supseteq W$ ($\beta v\alpha u=
\beta\alpha(\alpha^{-1}v\alpha)u\in\beta\alpha U$ for any $u\in U$ and 
$v\in V$, so the composition is well-defined). We endow ${\mathfrak T}$ 
with the maximal topology, i.e. we assume that any sieve is covering. 
Then the sheaves of sets, groups, etc. on ${\mathfrak T}$ are identified 
with the smooth $H$-sets, groups, etc.: 
${\mathcal F}\mapsto\lim\limits_{\overrightarrow{U\in B}}
{\mathcal F}(U)$ (this is a smooth $H$-set, since its arbitrary element 
belongs to the image of some ${\mathcal F}(U)$ and the $U$-action 
on it is trivial by definition) and $W\mapsto(U\mapsto W^U)$. 

{\sc Examples.} 1. If $H$ is discrete and $B=\{1\}$ then there is exactly one 
object $\ast$ in ${\mathfrak T}$ and ${\rm Hom}_{{\mathfrak T}}(\ast,\ast)=H$. 

2. In notation of Introduction, the set ${\rm Hom}_{{\mathfrak T}(G,B)}
(G_{F|L},G_{F|K})=\{h\in G~|~h(L)\subseteq K\}/G_{F|L}$ consists of all 
field embeddings $L\stackrel{/k}{\hookrightarrow}K$ over $k$. 
(Here $B$ is the set of open subgroups 
in $G$ of type $G_{F|L}$, where $L|k$ is an extension of finite type.) 

\subsection{Basic structures on ${\mathcal S}m_H(E)$} \label{prqmye-proizved} 
There are direct sums, direct products and the inner ${\mathcal H}om$ functor 
in the categories ${\mathcal S}m_H(E)$ and ${\mathcal I}_G$. They are the 
smooth parts of the corresponding functors on the category of $E$-vector 
spaces, and therefore, these functors on ${\mathcal I}_G$ coincide with the 
restrictions of the corresponding functors on ${\mathcal S}m_G({\mathbb Q})$. 
Namely, the direct sums are direct sums in the category of abelian 
groups; the direct product of a family of objects is the smooth part of its 
set-theoretic direct product: $\prod_{\alpha}^{{\mathcal S}m_H}W_{\alpha}:=
\bigcup_{U\in B}\prod_{\alpha}W_{\alpha}^U$, and ${\mathcal H}om(W_1,W_2)
:=\lim\limits_{\overrightarrow{U\in B}}{\rm Hom}_{E[U]}(W_1,W_2)$. 

The usual tensor structure on the category of $E$-vector spaces induces 
a tensor structure on the category ${\mathcal S}m_H(E)$. The functor 
${\mathcal H}om(W,-)$ is right adjoint to the functor $-\otimes_EW$: 
\begin{multline*} {\rm Hom}_{{\mathcal S}m_H(E)}(W_1\otimes_EW,W_2)
={\rm Hom}_{E[H]}(W_1\otimes_EW,W_2)={\rm Hom}_{E[H]}(W_1,
{\rm Hom}_E(W,W_2))=\\ ={\rm Hom}_{E[H]}(W_1,{\mathcal H}om(W,W_2))
={\rm Hom}_{{\mathcal S}m_H(E)}(W_1,{\mathcal H}om(W,W_2))\end{multline*}  
for any $W_1,W_2,W\in{\mathcal S}m_H(E)$. The tensor product on 
${\mathcal S}m_G({\mathbb Q})$ does not preserve ${\mathcal I}_G$, 
but it has a unique modification $\otimes_{{\mathcal I}}$ with the 
above adjunction property, cf. \cite[\S6.4]{repr}. The associativity 
of $\otimes_{{\mathcal I}}$, which is equivalent to the identity 
${\mathcal H}om(W_1\otimes W_2,-)=
{\mathcal H}om(W_1\otimes_{{\mathcal I}}W_2,-)$ on ${\mathcal I}_G$, 
is not yet known. 

\subsection{Generalities on Hecke algebras} \label{Hecke-algebras}
Define ${\mathbb D}_E(H):=
\!\lim\limits_{\phantom{_U}\longleftarrow_U}\!E[H/U]$, where the 
projective system is formed with respect to the projection $E[H/V]
\stackrel{r_{VU}}{\longrightarrow}E[H/U]$ and $H/V\longrightarrow H/U$, 
induced by the inclusions $V\subset U$ of open subgroups of $H$. In other 
words, ${\mathbb D}_E(H):=E\widehat{\otimes}{\mathbb D}_{{\mathbb Q}}(H)$. 

Any element $\nu\in{\mathbb D}_E(H)$ can be considered as an $E$-valued 
`oscillating' measure on $H$ (for which all open subgroups and 
their translates are measurable). In particular, for any $\sigma\in H$ and 
an open subgroup $U$ let the value $\nu(\sigma U)$ of the measure $\nu$ 
on the set $\sigma U$ be the $[\sigma U]$-coefficient 
in the image of $\nu$ in $E[H/U]$. 

For each smooth semilinear $E$-representation $W$ of $H$ define a pairing 
${\mathbb D}_E(H)\times W\longrightarrow W$ by $(\nu,w)\longmapsto
\sum_{\sigma\in H/U}\nu(\sigma U)\cdot\sigma w$, where $U$ is an arbitrary 
open subgroup in the stabilizer of $w$, e.g., $U={\rm Stab}_w$. 
Clearly, the result is independent of the choice of $U$. This determines 
a ${\mathbb D}_E(H)$-module structure on $W$. When $W=E[H/U]$, this pairing 
is compatible with the projections $r_{VU}$, so it gives rise to a pairing 
${\mathbb D}_E(H)\times\lim\limits_{\phantom{_U}\longleftarrow_U}\!E[H/U]
\longrightarrow\lim\limits_{\phantom{_U}\longleftarrow_U}\!E[H/U]=
{\mathbb D}_E(H)$, and thus, to an associative multiplication ${\mathbb D}_E(H)
\times{\mathbb D}_E(H)\stackrel{\ast}{\longrightarrow}{\mathbb D}_E(H)$, 
extending the convolution of the compactly supported measures. 
\label{support}
(The {\sl support} of $\nu$ is the minimal closed subset $S$ in the 
semi-group $\widehat{H}:=\lim\limits_{\phantom{_U}\longleftarrow_U}H/U$ 
such that $\nu(\sigma U)=0$ for any $\sigma U$ that does not meet $S$.) 

Clearly, the construction of ${\mathbb D}_E(H)$ is functorial 
in the sense that the following data \begin{itemize} 
\item a totally disconnected group $H'$, 
\item a field $E'$ endowed with a smooth $H'$-action, 
\item a continuous homomorphism $\varphi:H\longrightarrow H'$, 
\item a smooth 1-cocycle $\theta:H\longrightarrow(E')^{\times}$, 
i.e., $\theta|_U=1$ for some open subgroup $U\subset H$ and 
$\theta(hh')=\theta(h)\varphi(h)(\theta(h'))$, 
\item a field embedding $\lambda:E\longrightarrow E'$, 
which is compatible with the $H$-action, \end{itemize} induces a 
continuous homomorphism of algebras ${\mathbb D}_E(H)\longrightarrow
{\mathbb D}_{E'}(H')$, $a[h]\mapsto\lambda a\cdot\theta(h)[\varphi h]$. 

\label{def-Hecke}
{\sl The Hecke algebra} of a pair $(H,U)$, where $U$ is 
a compact subgroup of $H$, is the subalgebra 
${\mathcal H}_E(H,U):=h_U\ast{\mathbb D}_E(H)\ast h_U$ in 
${\mathbb D}_E(H)$ of $U$-biinvariant measures. Here $h_U$ is the Haar 
measure on $U$, defined by the system $(h_U)_V=[U:U\bigcap V]^{-1}
\sum_{\sigma\in U/U\bigcap V}[\sigma V]\in{\mathbb Q}[H/V]$ for all 
open subgroups $V\subset H$. The element $h_U$ is the unity of the 
algebras ${\mathcal H}_{{\mathbb Q}}(H,U)\subseteq{\mathcal H}_E(H,U)$ 
and $h_Uh_{U'}=h_U$ for any closed subgroup $U'\subseteq U$. 
For each smooth $E$-representation $W$ of $H$ the Hecke algebra 
${\mathcal H}_E(H,U)$ acts on the space $W^U$, since $W^U=h_U(W)$. 

This definition of the Hecke algebra is equivalent to the standard one 
in the case of locally compact $H$, open compact $U$ and a 
characteristic zero field $E$ endowed with the trivial $H$-action, 
and ${\mathcal H}_E(H,U)$ acts on the space $W^U$ in the usual way 
for any smooth $E$-representation $W$ of $H$, cf. \cite{bz}. 

If $H$ is the automorphism group of an algebraically closed extension of $k$ 
of finite transcendence degree $n$ then the Hecke algebras become the algebras 
of non-degenerate correspondences on some $n$-dimensional $k$-varieties, 
cf. Appendix \ref{ind-comp}. 

{\sc Remarks.} \label{semi-simplicity-general}
1. There is the following evident modification to arbitrary 
totally disconnected groups of the standard semi-simplicity or 
irreducibility criteria, cf. \cite[Proposition 2.10]{bz}. 

Let $T$ be a filtering family of compact subgroups of $H$, i.e. such 
that any open subgroup contains an element of $T$, e.g. the set of open 
subgroups of a given compact subgroup of $H$. Then a smooth 
$E$-representation $W$ of $H$ is irreducible, resp. semi-simple, if 
and only if the ${\mathcal H}_E(H,U)$-module $W^U$ is irreducible, 
resp. semi-simple, for each compact subgroup $U\in T$. 

2. One uses the centres of the Hecke algebras to decompose the category of 
smooth representations. However, in the case of $H=G$ one can show (in a way 
similar to that of \cite[Appendix A, Theorem A.4]{repr}) that for any compact 
subgroup $K$ in $G$ the centre of the Hecke algebra ${\mathcal H}_E(K)$ of 
the pair $(G,K)$ coincides with $E\cdot h_K$, i.e., consists of scalars. 

\subsection{Pull-back functors} 
If $\varphi:H_2\longrightarrow H_1$ is a homomorphism 
with a dense image then the pull-back functor $\varphi^{-1}:
{\mathcal S}m_{H_1}\longrightarrow H_2$-mod is fully faithful. 
({\it Proof.} Let $W_1,W_2\in{\mathcal S}m_{H_1}$, 
$\alpha\in{\rm Hom}_{H_2}(\varphi^{-1}W_1,\varphi^{-1}W_2)$, 
$v\in W_1$ and $\sigma\in H_1$. Let $S$ be the common stabilizer of the 
elements $v$ and $\alpha(v)$. Choose some element 
$\sigma'\in\varphi^{-1}(\sigma S)\subseteq H_2$. Then 
$\alpha(\sigma v)=\alpha(\sigma'v)=\sigma'\alpha(v)=\sigma\alpha(v)$. \qed) 

If $\varphi$ is continuous then $\varphi^{-1}$ factors through 
$\varphi^{\ast}:{\mathcal S}m_{H_1}\longrightarrow{\mathcal S}m_{H_2}$. 

If the homomorphism $\varphi$ is continuous and with dense image then 
the functor $\varphi^{\ast}$ admits a right adjoint $\varphi_{\ast}:
W\mapsto\bigcup_UW^{U\times_{H_1}H_2}$, where $U$ runs over open subgroups 
of $H_1$. In particular, $\varphi^{\ast}$ preserves the irreducibility. 
({\it Proof.} The $H_1$-action on $\varphi_{\ast}W$ is defined as follows. 
If $w\in W^{\varphi^{-1}(U)}$ and $\sigma\in H_1$ then 
$\sigma w:=\sigma'w$, where $\sigma'\in H_2$ and 
$\varphi(\sigma')\in\sigma U$, which is independent of $\sigma'$. \qed) 

{\sc Example.} (A `dense' locally compact `subgroup' ${\mathfrak G}$ 
of $G$.) Let $\{x_1,x_2,\dots\}$ be a transcendence base of $F|k$. 
Set $L_m:=k(x_m,x_{m+1},\dots)\subset F$. Then $L_{\bullet}=(L_1\supset 
L_2\supset L_3\supset\dots)$ is a descending sequence of subfields in $F$. 
Set ${\mathfrak G}={\mathfrak G}_{L_{\bullet}}:=\bigcup\limits_{m\ge 1}
G_{F|L_m}$. We take the set $\{G_{F|LL_1}\}$ of subgroups for all 
subfields $L$ in $F|k$ of finite type as a base of open subgroups. 

Geometrically (in a sense, analogous to the dominant topology), 
this corresponds to an inverse system of infinite-dimensional irreducible 
$k$-varieties given by finite systems of equations. They are related by 
dominant morphisms affecting only finitely many coordinates. 

Then \begin{itemize}\item ${\mathfrak G}$ is locally compact (since 
$F$ is algebraic over $L_1$), but is not unimodular; 
\item the inclusion ${\mathfrak G}$ into $G$ is continuous with 
dense image (since $\bigcap_{m\ge 1}\overline{L_m}=k$). \end{itemize}
More details can be found in \cite{rms}. 

The forgetful functor ${\mathcal S}m_G\longrightarrow
{\mathfrak G}\mbox{-}{\rm mod}$ is fully faithful, preserves the 
irreducibility, factors through $r:{\mathcal S}m_G\longrightarrow
{\mathcal S}m_{{\mathfrak G}}$, and $r$ admits a right adjoint: 
$W\mapsto\bigcup_L\bigcap_{m\ge 1}W^{G_{F|LL_m}}$, where $L$ 
runs over the set of all subfields of finite type in $F|k$. 

\subsection{Count of cyclic and irreducible objects of ${\mathcal S}m_H(E)$} 
\label{count-cyclic-irr}
Recall (e.g., \cite[Ch. III, \S1]{milne}) that a family $(A_j)_{j\in J}$ 
of objects in a category is {\sl generating} if for any injection 
$\alpha:A'\hookrightarrow A$ which is not an isomorphism there exist 
an index $j\in J$ and a morphism $A_j\longrightarrow A$ that 
does not factor through $\alpha$. 

The objects $E[H/U]$, where $U\in B$, form a generating system 
of ${\mathcal S}m_H(E)$, i.e., any smooth cyclic $E$-semilinear 
representation of $H$ is a quotient of $E[H/U]$ for some $U\in B$. 

One can modify the $H$-module structure on $E[H/U]$ by the rule 
$\sigma:[h]\mapsto f(\sigma,h)\cdot[\sigma h]$ for some function 
$f:H\times H/U\longrightarrow E^{\times}$. The associativity 
constraint gives the condition 
$f(\tau\sigma,h)=\tau f(\sigma,h)\cdot f(\tau,\sigma h)$. In other 
words, $g(\sigma\tau,h)=g(\sigma,h)\cdot\sigma g(\tau,\sigma^{-1}h)$, 
where $g(\sigma,h)=\sigma f(\sigma^{-1},h)$, i.e., 
$g\in Z^1(H,{\rm Maps}(H/U,E^{\times}))$ is a 1-cocycle on $H$. 
The smoothness constraint amounts to the condition that for any 
$h\in H/U$ there exists an open subgroup $V\subset H$ such that 
$f(\sigma,h)=1$ (equivalently, $g(\sigma,h)=1$) for any $\sigma\in V$. 
We denote by $E[H/U]^{(g)}$ the resulting object of ${\mathcal S}m_H(E)$. 

Let $\delta:{\rm Maps}(H/U,E^{\times})\longrightarrow 
Z^1(H,{\rm Maps}(H/U,E^{\times}))$, $\varphi\longmapsto
[(\sigma,h)\mapsto\varphi(h)^{-1}\cdot\sigma\varphi(\sigma^{-1}h)]$, 
be the coboundary homomorphism. Any coboundary is smooth: 
$(\sigma,h)\mapsto\varphi(h)^{-1}\varphi(h)=1$ for any 
$\sigma\in{\rm Stab}_{\varphi(h)}\cap hUh^{-1}$. If two 1-cocycles 
are cohomological to each other then the corresponding objects of 
${\mathcal S}m_H(E)$ are isomorphic: $E[H/U]^{(g)}
\stackrel{\sim}{\longrightarrow}E[H/U]^{(g\cdot\delta\varphi)}$, 
$[h]\mapsto\varphi(h)\cdot[h]$. 

{\sc Example.} Suppose that $\varphi$ lifts to a 1-cocycle 
$H\longrightarrow E^{\times}$, and thus, corresponds to a one-dimensional 
object ${\mathcal L}$ of ${\mathcal S}m_H(E)$. Then 
$(\delta\varphi)(\sigma,h)=\varphi(h)^{-1}\cdot\sigma\varphi(\sigma^{-1}h)
=\varphi(h)^{-1}\varphi(\sigma)^{-1}\varphi(h)=\varphi(h)^{-1}$, 
so $f(\sigma,h)=\sigma\varphi(\sigma^{-1})^{-1}=\varphi(\sigma)$, 
and therefore, $E[H/U]^{(g\cdot\delta\varphi)}\cong 
E[H/U]^{(g)}\otimes_E{\mathcal L}$. 

As the representations $E[H/U]^{(g)}$ and $E[H/U]^{(g\cdot\delta\varphi)}$ 
are isomorphic, this implies that for any irreducible $W\in{\mathcal S}m_H(E)$ 
and any one-dimensional ${\mathcal L}\in{\mathcal S}m_H(E)$ with 
${\mathcal L}^U\neq 0$ the multiplicities of $W$ and of 
$W\otimes_E{\mathcal L}$ in $E[H/U]^{(g)}$ coincide. 

It is likely, however, that these multiplicities are infinite, 
when $H$ is the automorphism group of a non-trivial algebraically 
closed extension of $k$ of finite transcendence degree, cf. 
\cite[Remark on p.217]{repr} (or \cite[p.1162]{rms}). 
More remarks on representations $E[G/U]$ are in Appendix \ref{ind-comp}. 

\begin{proposition} \label{neschet-nepriv} \begin{enumerate}
\item \label{obx-verh-cikl} There are at most 
$\max(|B|,\sup\limits_{U\in B}2^{\max(|H/U|,|E|)})$ isomorphism 
classes of smooth cyclic $E$-semilinear representations of $H$. 
\item \label{nizhn-nepriv-slab} Suppose that $U$ is an open subgroup of $H$, 
and there are $>\sup\limits_{U'\in B}\max(|H/U'|,|E|)$ isomorphism classes 
of irreducible objects of ${\mathcal S}m_U(E)$. Then the group $H$ admits  
at least as many isomorphism classes of smooth irreducible representations 
as $U$ does. 
\item \label{nizhn-nepriv-siln} Let $U$ be an open subgroup of $H$. Suppose 
that there are at least $\max(|B|,\sup\limits_{U'\in B}2^{\max(|H/U'|,|E|)})$ 
isomorphism classes of irreducible objects of ${\mathcal S}m_U(E)$. 
Then the cardinality of the set of isomorphism classes of irreducible 
objects of ${\mathcal S}m_H(E)$ is equal to 
$\max(|B|,\sup\limits_{U'\in B}2^{\max(|H/U'|,|E|)})$. 
\end{enumerate}
\end{proposition}
{\it Proof.} As $E[S]$ is dominated by $\bigcup_{N\ge 1}E^N\times S^N$, 
$(a_1,\dots,a_N;s_1,\dots,s_N)\mapsto\sum_{i=1}^Na_i[s_i]$, one has 
$|E[S]|=\max(|S|,|E|)$, if $E$ and $S$ are infinite. Then for any 
subgroup $U$ of $H$ there are $\le 2^{\max(|H/U|,|E|)}$ quotients 
of the representation $E[H/U]$, which proves (\ref{obx-verh-cikl}).

Let $U$ be an open subgroup of $H$, and $\varphi$ be a smooth cyclic 
$E$-representation of $U$. Let $W$ be any quotient of the cyclic 
representation $E[H]\otimes_{E[U]}\varphi$ of $H$. Then there are 
$\le|W|\le\sup_{U'\in B}|E[H/U']|\le\sup_{U'\in B}\max(|H/U'|,|E|)$ 
isomorphism classes of cyclic subrepresentations of $U$ in $W$, one 
of which is $\varphi$, if $\varphi$ is irreducible. For each smooth 
irreducible representation $\varphi$ of $U$ choose an irreducible 
quotient $W_{\varphi}$ of the representation $E[H]\otimes_{E[U]}\varphi$. 

We say that smooth irreducible representations of $U$ are equivalent 
(notation: $\varphi\sim\psi$) if $W_{\varphi}\cong W_{\psi}$. As 
$|W_{\varphi}|\le\sup_{U'\in B}\max(|H/U'|,|E|)$, the cardinalities of the 
equivalence classes do not exceed $\sup\limits_{U'\in B}\max(|H/U'|,|E|)$. 
Therefore, $|\{\varphi\}/\sim|=|\{\varphi\}|$ under assumption of 
(\ref{nizhn-nepriv-slab}). 

Under assumption of (\ref{nizhn-nepriv-siln}), the set of 
isomorphism classes of smooth irreducible representations of $U$ and 
the set of equivalence classes of smooth irreducible representations 
of $U$ have the same cardinality, and the lower bound 
$\ge\max(|B|,\sup\limits_{U'\in B}2^{\max(|H/U'|,|E|)})$ of this 
cardinality coincides with the upper bound of the cardinality of 
the set of isomorphism classes of smooth irreducible representations 
of $H$ from (\ref{obx-verh-cikl}), so the group $H$ admits 
precisely $\max(|B|,\sup\limits_{U'\in B}2^{\max(|H/U'|,|E|)})$ 
isomorphism classes of smooth irreducible representations. \qed 

\vspace{4mm}

{\sc Examples.} 1. Let $H$ be the automorphism group of a non-trivial 
algebraically closed field extension $F$ of $k$ of transcendence degree 
at most $|k|$, e.g., countable. Let $B$ be the set of stabilizers of 
finite subsets in $F$. Then $|B|=\sup\limits_{U'\in B}|H/U'|=|k|$, 
and thus, according to Proposition \ref{neschet-nepriv} 
(\ref{obx-verh-cikl}), there are $\le 2^{\max(|k|,|E|)}$ isomorphism 
classes of cyclic objects of ${\mathcal S}m_H(E)$. Let $U$ be 
the (setwise) stabilizer in $H$ of an algebraically closed 
subfield $F'\neq k$ in $F|k$ of a finite transcendence degree. 

Note, that if the $H$-action on $E$ is trivial then there are 
\begin{itemize} \item 
$|{\rm Hom}({\mathbb Q}^{\times}_+,E^{\times})|=\max(2^{|{\mathbb N}|},|E|)$ 
one-dimensional representations of $U$, that factor through the composition 
of the surjective continuous homomorphisms $U\longrightarrow G_{F'|k}
\stackrel{\chi}{\longrightarrow}{\mathbb Q}^{\times}_+$, where $\chi$ 
is the modulus of $G_{F'|k}$; 
\item $\ge|k|$ smooth irreducible $E$-representations of $H$: for each 
elliptic curve $A$ over $k$ without complex multiplication the smooth 
representation $(A(F)/A(k))\otimes E$ of $H$ is irreducible. 

\end{itemize}

By Proposition \ref{neschet-nepriv} (\ref{nizhn-nepriv-siln}), there are 
exactly $2^{|{\mathbb N}|}$ smooth irreducible $E$-representations of $G$, 
if $k$ and $E$ are countable.

Proposition \ref{neschet-nepriv} shows, there are `too many' 
($\ge\max(2^{|{\mathbb N}|},|k|,|E|)$) smooth irreducible 
representations of $G$. (This is one of the reasons to study rather 
${\mathcal I}_G$  than ${\mathcal S}m_G(E)$, where the objects are 
supposed to be more controllable, since it is expected 
that they are of `cohomological nature'.)

2. If $H$ is locally compact, but not unimodular, then there are 
$\ge\max(2^{{\rm rk}(\chi_H(H))},|E|)$ irreducible representations. 
E.g., if $H$ is the automorphism group of a non-trivial algebraically 
closed field extension of $k$ of finite transcendence 
degree, we get the bound $\ge\max(2^{|{\mathbb N}|},|E|)$.

\subsection{Injectives in ${\mathcal S}m_H(E)$ and in ${\mathcal I}_G$}
For any $E$-vector space $V$, let ${\rm Maps}(H,V)$ be the module of 
$V$-valued functions on $H$. We endow ${\rm Maps}(H,V)$ with the following 
$E$-vector space structure: $(\lambda f)(x)=x\lambda\cdot f(x)$ for any 
$\lambda\in E$. The $H$-module structure on ${\rm Maps}(H,V)$ is given by the 
right translations of the argument $h:f(x)\mapsto f(xh)$ for any $h\in H$. 

Define $I(V)$ as the smooth part of ${\rm Maps}(H,V)$, i.e., 
$I(V):=\lim\limits_{\overrightarrow{U'\in B}}{\rm Maps}(H/U',V)$. 
Clearly, $I(V)$ is a $E$-vector subspace of ${\rm Maps}(H,V)$ invariant 
under the $H$-action. More generally, for any subgroup $U\subset H$ 
denote by ${\rm Ind}^H_U:{\mathcal S}m_U(E)\longrightarrow
{\mathcal S}m_H(E)$ the smooth induction functor $V\longmapsto
\lim\limits_{\overrightarrow{U'\in B}}{\rm Maps}_U(H/U',V)$, 
where ${\rm Maps}_U(H/U',V)$ is the group of $U$-equivariant maps. 

In fact, ${\rm Ind}^H_U$ is right adjoint to the forgetful functor 
${\rm Res}^H_U:{\mathcal S}m_H(E)\longrightarrow{\mathcal S}m_U(E)$: 
\begin{equation}\label{adjointness} {\rm Hom}_{{\mathcal S}m_H(E)}
(A,{\rm Ind}^H_U(V))={\rm Hom}_{E[H]}(A,{\rm Maps}_U(H,V))={\rm Hom}
_{{\mathcal S}m_U(E)}({\rm Res}^H_UA,V)\end{equation} for any 
$A\in{\mathcal S}m_H(E)$ and $V\in{\mathcal S}m_U(E)$, where 
${\rm Hom}_{E[H]}(A,{\rm Maps}_U(H,V))\ni\varphi\mapsto
(a\mapsto\varphi(a)(1))\in{\rm Hom}_{E[U]}(A,V)$; 
${\rm Hom}_{E[U]}(A,V)\ni\psi\mapsto(a\mapsto(h\mapsto\psi(ha)))
\in{\rm Hom}_{E[H]}(A,{\rm Maps}_U(H,V))$.  

Therefore, if $V$ is an injective object in ${\mathcal S}m_U(E)$ then 
${\rm Ind}^H_U(V)$ is an injective object in ${\mathcal S}m_H(E)$. 
In particular, it follows from the semi-simplicity of the category 
of $E$-vector spaces that its arbitrary object is injective, and thus, 
$I(V)={\rm Ind}^H_{\{1\}}(V)$ is also injective. If, moreover, $V$ is a 
smooth $H$-module then $V\longrightarrow I(V)$, $v\mapsto(h\mapsto hv)$, 
is a monomorphism in ${\mathcal S}m_H(E)$. This means that there are 
enough injectives in ${\mathcal S}m_H(E)$. 

\begin{lemma}[\cite{milne} Ch.III, Lemma 1.3] If in an abelian category 
there are direct sums, products, a generating family and all filtered 
direct limits are exact then there are enough injectives. \end{lemma}

For a subgroup $U$ of $H$ denote by $H^q_{{\mathcal S}m_H}(U,-)$ 
the $q$-th derived functor\footnote{As there are enough injectives in 
${\mathcal S}m_H$ the derived functors of left exact functors are defined.} 
of $H^0(U,-)$ on ${\mathcal S}m_H$. Evidently, $H^q_{{\mathcal S}m_H}(U,-)
=H^q_{{\mathcal S}m_H}(\overline{U},-)$, where $\overline{U}$ 
is the closure of $U$ in $H$. 

In other words, $H^j_{{\mathcal S}m_H}(H,-)=
{\rm Ext}^j_{{\mathcal S}m_H}({\mathbb Q},-)$ 
can be defined using the smooth cochains. 

\begin{lemma} \label{open-support} There is a morphism of 
functors $H^{\ast}_{{\mathcal S}m_H}(U,-)\longrightarrow 
H^{\ast}_{{\mathcal S}m_U}(U,-)$ on ${\mathcal S}m_H$, 
which is an isomorphism if $U$ is open. \end{lemma} 
{\it Proof.} As the right ${\mathbb Z}[U]$-module 
${\mathbb Z}[H]$ is free, the co-induction functor ${\mathbb Z}[H]
\otimes_{{\mathbb Z}[U]}-:{\mathcal S}m_U\longrightarrow
{\mathcal S}m_H$ is exact. On the other hand, it is left adjoint 
to the forgetful functor ${\rm Res}^H_U:{\mathcal S}m_H
\longrightarrow{\mathcal S}m_U$, and therefore, for any injective 
object $I\in{\mathcal S}m_H$ the functor 
${\rm Hom}_U(-,{\rm Res}^H_U(I))={\rm Hom}_H({\mathbb Z}[H]
\otimes_{{\mathbb Z}[U]}-,I)$ on ${\mathcal S}m_U^{{\rm op}}$ 
is exact. In other words, the forgetful functor ${\rm Res}^H_U$ 
transforms injectives of ${\mathcal S}m_H$ to injectives of 
${\mathcal S}m_U$, and thus, an injective resolution of $W$ 
in ${\mathcal S}m_H$ becomes an injective resolution of 
${\rm Res}^H_U(W)$ in ${\mathcal S}m_U$. \qed

\vspace{4mm}

{\sc Remark.} Evidently, the morphism of 
functors $H^{\ast}_{{\mathcal S}m_H}(U,-)\longrightarrow 
H^{\ast}_{{\mathcal S}m_U}(U,-)$ on ${\mathcal S}m_H$, 
is not an isomorphism for arbitrary $U$: if $H$ is topologically simple, 
e.g. $H=G$, and $U$ is locally compact and not unimodular, e.g. $U=G_{F|F'}$ 
for a subfield $F'\subset F$, over which $F$ is of finite transcendence 
degree, then $H^{\ast}_{{\mathcal S}m_H}(U,{\mathbb Q})=0\neq 
H^{\ast}_{{\mathcal S}m_U}(U,{\mathbb Q})$. 

\begin{lemma} \label{limit-support} Let $U\subset H$ be a subgroup, 
and $\{U\subset U_{\alpha}\subset H\}_{\alpha}$ be a partially ordered 
collection of subgroups, which is filtering for the neighbourhoods of $[1]
\in H/\overline{U}$, i.e. any neighbourhood of $[1]\in H/\overline{U}$ 
contains $U_{\alpha}/U_{\alpha}\cap\overline{U}\subset H/\overline{U}$ 
for some $\alpha$. Then 
$\lim\limits_{_{\alpha}\longrightarrow\phantom{_{\alpha}}}
H^q_{{\mathcal S}m_H}(U_{\alpha},V)\stackrel{\sim}{\longrightarrow}
H^q_{{\mathcal S}m_H}(U,V)$. \end{lemma} 
{\it Proof.} Let $I^{\bullet}(V)$ be an injective resolution of $V$. 
By definition, $H^q_{{\mathcal S}m_H}(U,V)=H^q\left(I^{\bullet}(V)^U\right)$. 
As all $I^i(V)$ are smooth, $I^{\bullet}(V)^U=
\lim\limits_{_{\alpha}\longrightarrow\phantom{_{\alpha}}}
I^{\bullet}(V)^{U_{\alpha}}$. As the direct limits commute with the 
cohomology, $H^q\left(\lim\limits_{_{\alpha}\longrightarrow\phantom{_{\alpha}}}
I^{\bullet}(V)^{U_{\alpha}}\right)=
\lim\limits_{_{\alpha}\longrightarrow\phantom{_{\alpha}}}
H^q\left(I^{\bullet}(V)^{U_{\alpha}}\right)=
\lim\limits_{_{\alpha}\longrightarrow\phantom{_{\alpha}}}
H^q_{{\mathcal S}m_H}(U_{\alpha},W)$, which completes the proof. \qed

\vspace{4mm}

{\sc Example.} For any subgroup $U$ of $H$ the collection $\{\langle U,
U'\rangle\}_{U'\in B}$ of subgroups of $H$, generated by $U$ and the 
elements of $B$, satisfies the assumptions of Lemma \ref{limit-support}. 

In this case $\lim\limits_{_{\alpha}\longrightarrow\phantom{_{\alpha}}}
H^q_{{\mathcal S}m_{U_{\alpha}}}(U_{\alpha},W)\stackrel{\sim}
{\longrightarrow}H^q_{{\mathcal S}m_H}(U,W)$ by Lemma \ref{open-support}. 

\begin{proposition} For any $j\ge 0$ the restriction to ${\mathcal S}m_H(E)$ 
of the functor ${\rm Ext}^j_{{\mathcal S}m_H}({\mathbb Q},-)$ coincides with 
${\rm Ext}^j_{{\mathcal S}m_H(E)}(E,-)$. Clearly, $H^0(H,-)={\rm Hom}_H
({\mathbb Q},-)={\rm Hom}_{{\mathcal S}m_H(E)}(E,-)$ on ${\mathcal S}m_H(E)$. 
\end{proposition}
{\it Proof.} This is a particular case of 
\cite[Expos\'{e} V, Corollaire 3.5]{sga-4-2} and follows from the fact 
that the objects $W$ of ${\mathcal S}m_H(E)$ admit canonical injective 
resolutions $I(W)\to I(I(W)/W)\to I(I(I(W)/W)/I(W))\to\dots $ such that 
the forgetful functor ${\mathcal S}m_H(E)\longrightarrow{\mathcal S}m_H$ 
transforms them to injective resolutions in ${\mathcal S}m_H$. 
The injectivity follows from the above adjointness property 
(\ref{adjointness}). \qed 

\begin{corollary} There are enough injectives in 
${\mathcal S}m_H(E)$ and in ${\mathcal I}_G$. \end{corollary} 
{\it Proof.} Direct sums and direct products in ${\mathcal S}m_H(E)$ and in 
${\mathcal I}_G$ were already defined in Appendix \ref{prqmye-proizved}. 
Direct limits in all these categories are direct limits in the 
category of abelian groups, and therefore, they are exact. 
For ${\mathcal S}m_H(E)$ the family $(A_U=E[H/U])_{U\in B}$ is 
generating. (Clearly, the direct sum or direct product of all objects in 
a generating family of ${\mathcal S}m_H(E)$ is a generator.) In the case 
$E=F$ for each positive integer $m$ fix a field extension $K_m$ of $k$ of a 
finite transcendence degree $\ge m$. Then $(A_m=F[G/G_{F|K_m}])_{m\ge 1}$ 
is another generating family for ${\mathcal S}m_G(F)$. The functor 
${\mathcal I}$, left adjoint to the inclusion 
${\mathcal I}_G\hookrightarrow{\mathcal S}m_G$, cf. \cite{repr}, 
transforms any generating family for 
${\mathcal S}m_G$ to a generating family for ${\mathcal I}_G$. \qed 

\subsection{Projectives in ${\mathcal S}m_H(E)$} \label{no-proj} 
Recall, cf. p.\pageref{support}, that 
$\widehat{H}:=\lim\limits_{\phantom{_U}\longleftarrow_U}H/U$ is a semi-group. 
\begin{lemma} \label{prakticheski-loc-comp} 
The following conditions are equivalent: \begin{enumerate} 
\item \label{equiv-loc-comp} there is a locally compact group $H'$ 
and an equivalence of categories ${\mathcal S}m_H(E)
\stackrel{\sim}{\longrightarrow}{\mathcal S}m_{H'}(E)$; 
\item \label{essent-loc-comp} $\widehat{H}$ is a locally compact group, 
so $\widehat{H}=\lim\limits_{\phantom{_U}\longleftarrow_U}U\backslash H$, 
and $H\longrightarrow\widehat{H}$ is an embedding with dense image such 
that the topology of $H$ is induced by the topology of $\widehat{H}$; 
\item \label{exists-subgroup-fin-ind} there exists an open 
subgroup of $H$ such that its any open subgroup is of finite index; 
\item \label{exists-proj} there exists a non-zero projective object 
in the category ${\mathcal S}m_H(E)$ for some $E$; 
\item \label{exists-enough-proj} there exist enough projective objects in 
the category ${\mathcal S}m_H(E)$ for any $E$. \end{enumerate} \end{lemma} 
{\it Proof.} It is evident that (\ref{essent-loc-comp}) and 
(\ref{exists-subgroup-fin-ind}) are equivalent. The implication 
(\ref{exists-enough-proj})$\Rightarrow$(\ref{exists-proj}) is trivial. 

(\ref{exists-subgroup-fin-ind})$\Rightarrow$(\ref{exists-enough-proj}). 
Let $U$ be an open subgroup of $H$ such that its any open subgroup is of 
finite index. Then the objects $E[H/U']$ are projective for all open 
subgroups $U'$ of $U$ and form a system of generators of ${\mathcal S}m_G(E)$. 

(\ref{exists-proj})$\Rightarrow$(\ref{exists-subgroup-fin-ind}). 
Let $W\in{\mathcal S}m_H(E)$ be a projective object. 
Choose a generating system $\{e_j\}_{j\in J}$ of the representation $W$. 
This gives rise to a surjective homomorphism 
$\bigoplus_{j\in J}E[H/{\rm Stab}_{e_j}]\stackrel{\pi}{\longrightarrow}W$. 
Fix an element $i_0\in J$ and for each $j\in J$ fix an open subgroup $U_j$ 
in ${\rm Stab}_{e_j}\cap{\rm Stab}_{e_{i_0}}$. As $W$ is projective, 
the composition of $\pi$ with the surjection $\bigoplus_{j\in J}
E[H/U_j]\longrightarrow\bigoplus_{j\in J}E[H/{\rm Stab}_{e_j}]$ splits, 
and therefore, there exists an element in $\bigoplus_{j\in J}E[H/U_j]$ 
with the same stabilizer as $e_{i_0}$. This implies that the space 
$E[H/U_j]^{{\rm Stab}_{e_{i_0}}}$ is non-zero for some $j$, and thus, 
the index of $U_j$ in ${\rm Stab}_{e_{i_0}}$ is finite. In other words, 
any open subgroup of ${\rm Stab}_{e_{i_0}}$ is of finite index. \qed 

\vspace{4mm}

{\sc Examples.} 1. In the case $H=G$ any open subgroup of $H$ contains 
an open subgroup of infinite index, and thus, there are no 
non-zero projective objects in the category ${\mathcal S}m_G(E)$. 

2. If $H$ is locally compact then there are enough projectives in 
${\mathcal S}m_H(E)$. 

\subsection{Left exact subfunctors of the forgetful functor from 
${\mathcal S}m_H(E)$ to the category of $E$-vector spaces} 
\label{leftexact-subf} 
The functors $H^0(U,-)$ on ${\mathcal S}m_H(E)$ for 
subgroups\footnote{Clearly, $H^0(U,-)=H^0(U',-)$, where $U'$ is the 
intersection of all open subgroups in $H$ containing $U$.} 
$U\subseteq H$ are examples of left exact subfunctors of the forgetful 
functor from ${\mathcal S}m_H(E)$ to the category of $E$-vector spaces. 

Fix a `sufficiently big' $E$-vector space $W_0$. If 
$\dim_EW_0\geqslant|H/U|$ for any open subgroup $U\subset H$ then any 
finitely generated representation of $H$ is embeddable into $I(W_0)$. 
\begin{proposition} There is a natural bijection between the $E$-vector 
subspaces in $I(W_0)$, invariant under the algebra 
${\rm End}_{E[H]}(I(W_0))={\rm Hom}_E(I(W_0),W_0)$, and left exact 
subfunctors of the forgetful functor from ${\mathcal S}m_H(E)$ to 
the category of $E$-vector spaces.\end{proposition} 
Namely, any functor $\Phi$ is sent to the subspace $\Phi(I(W_0))$. 
Conversely, as the $E$-vector space $W_0$ is `big enough', any 
`sufficiently small' object of ${\mathcal S}m_H(E)$ is embeddable into 
$I(W_0)$, and thus, one can define $\Phi(W):=W\bigcap\Phi(I(W_0))$, 
which is independent of the choice of the embedding. \qed 

\section{Non-degenerate generically finite correspondences and smooth 
$G$-modules coinduced from open subgroups} \label{ind-comp} 
Consider the category of smooth $k$-varieties with the morphisms, 
given by formal linear combinations of non-degenerate generically 
finite correspondences, i.e. irreducible subvarieties in the 
product of the source and the target, generically finite over a 
connected component of the source and dominant over a connected 
component of the target: ${\rm Hom}(X,Y)=Z^{\dim Y}(k(X)\otimes_kk(Y))$ 
for connected $X$ and $Y$. 

Then there is a full embedding of this category into the category of smooth 
representations of $G$, given by $X\mapsto Z^{\dim X}(k(X)\otimes_kF)=
{\mathbb Q}[\{k(X)\stackrel{/k}{\hookrightarrow}F\}]$. 
If we fix a $k$-field embedding of the function field $k(X)$ 
into $F$ then the module $Z^{\dim X}(k(X)\otimes_kF)$ of 
generic 0-cycles on $X_F$ becomes ${\mathbb Q}[G/G_{F|k(X)}]$. 
These $G$-modules are very complicated. 

\vspace{4mm}

One can extract dimension of $X$ out of $W=Z^{\dim X}(k(X)\otimes_kF)$: 
it is equal to the minimal $q\ge 0$ such that $W^{G_{F|\overline{L}}}\neq 0$, 
where ${\rm tr.deg}(L|k)=q$. Also, `birational motivic' invariants `modulo 
isogenies', such as ${\rm Alb}(X)$, $\Gamma(X,\Omega^{\bullet}_{X|k})
={\rm Hom}_G(C_{k(X)},\Omega^{\bullet}_{F|k})$, can be recovered 
from $W$, cf. \cite{repr,pgl}. However, we do not know, whether 
the birational type of $X$ is determined by $W$. 

The collection ${\rm JH}(X)$ of irreducible subquotients of the $G$-module 
of generic 0-cycles on $X$ over $F$ is a birational invariant of $X$. In 
this section we give examples of pairs of non-birational varieties $X$ and 
$Y$, with the same collections ${\rm JH}(X)$ and ${\rm JH}(Y)$. It is not 
even excluded that ${\rm JH}(X)$ depends only on $\dim X$ (and moreover, 
it is not even shown yet that it does depend on $X\neq{\rm Spec}(k)$). 

To show the inclusion ${\rm JH}(X)\supseteq{\rm JH}(Y)$, it suffices 
to construct a $G$-embedding of $Z^{\dim X}(k(Y)\otimes_kF)$ into a 
cartesian power of $Z^{\dim X}(k(X)\otimes_kF)$. For any generically 
finite map $X\longrightarrow Y$ the pull-back induces a desired embedding 
$Z^{\dim X}(k(Y)\otimes_kF)\hookrightarrow Z^{\dim X}(k(X)\otimes_kF)$. 
(This is a particular case of the following situation. Let 
$U\subset U'\subset H$ be subgroups, and let index of $U$ 
in $U'$ be finite. Then $[u]\mapsto\sum_{h\in H/U,~hU'=uU'}[h]$ 
gives a natural $H$-embedding $E[H/U']\hookrightarrow E[H/U]$.)

\begin{proposition}[\cite{repr}, Corollary 7.3] \label{example} 
Let $Z$ be a $k$-variety, $Z'$ be a generically twofold cover of $Z$, 
$X=Z\times{\mathbb P}^1$ and $Y=Z'\times{\mathbb P}^1$. Then there are 
$G$-embeddings 
$Z^{\dim X}(k(X)\otimes_kF)\hookrightarrow Z^{\dim X}(k(Y)\otimes_kF)$ and 
$Z^{\dim X}(k(Y)\otimes_kF)\hookrightarrow Z^{\dim X}(k(X)\otimes_kF)$, so 
${\rm JH}(X)={\rm JH}(Y)$. \end{proposition}
This results from the following combinatorial claim 
(\cite[Lemma 7.2]{repr}). 

{\it Let $H$ be a group and $U$ and $U'$ be subgroups of $H$ such that 
$U\bigcap U'$ is of index two in $U$: $U=(U\bigcap U')\bigcup\sigma
(U\bigcap U')$. Suppose that $\tau_1\cdots\tau_N\neq 1$ for any integer 
$N\ge 1$ and for any collection $\tau_1,\dots,\tau_N\in U'\sigma
\smallsetminus U$. Then the morphism of $E$-representations 
$E[H/U]\stackrel{[\xi]\mapsto[\xi\sigma]+[\xi]}
{-\!\!\!-\!\!\!-\!\!\!-\!\!\!-\!\!\!-\!\!\!\longrightarrow}E[H/U']$ 
of $H$ is injective.}

\begin{proposition}[\cite{repr}, 7.4] \label{example-2} Fix an 
odd integer $m\ge 1$, and let the affine $(m-1)$-dimensional $k$-variety 
$Y$ be given by equation $\sum_{j=1}^mx_j^d=1$, where $d\in\{m+1,m+2\}$. 
Let $X$ be the quotient of $Y$ by $\langle e_1e_2^2\cdots e_m^m\rangle$, 
where $e_ix_j=\zeta^{\delta_{ij}}\cdot x_j$ for a primitive $d$-th root 
of unity $\zeta$. Then ${\rm JH}(X)={\rm JH}({\mathbb P}^{m-1}_k)$. 
\end{proposition}

\begin{proposition} \label{example-ellipt-prqmaq} Let $U$ be a subgroup 
of $H$ and $g_1,\dots,g_N$ be involutions in $N_HU/U$, generating an 
infinite subgroup of $N_HU/U$. Then the natural map of $E$-representations 
$r:E[H/U]\longrightarrow\bigoplus_{j=1}^NE[H/\langle U,g_j\rangle]$ 
of $H$ is injective. 

In particular, if $g_1,\dots,g_N$ are rational involutions of a $k$-variety 
$X$, generating an infinite group, then the natural $G$-morphism 
$Z^{\dim X}(k(X)\otimes_kF)\longrightarrow\bigoplus_{j=1}^NZ^{\dim X}
(k(X)^{\langle g_j\rangle}\otimes_kF)$ is injective. \end{proposition}
{\it Proof.} If a non-zero 0-cycle $\alpha$ is in the kernel of $r$, 
and $P$ is a point in the support of $\alpha$, then the support of $\alpha$ 
contains the $\langle g_1,\dots,g_N\rangle$-orbit of the point $P$. 
As this orbit is infinite, but the support of $\alpha$ is finite, 
we get the contradiction, i.e., $\alpha=0$. \qed 

{\sc Examples.} 1. Let $X$ be an algebraic $k$-group, 
$g_1:x\mapsto x^{-1}$ and $g_2:x\mapsto h\cdot x^{-1}$, where $h\in X(k)$ 
is a point of infinite order. Then the $E$-representations $E[\{k(X)
\stackrel{/k}{\hookrightarrow}F\}]$ and $E[\{k(K(X))\stackrel{/k}
{\hookrightarrow}F\}]$ of $G$ have the same irreducible subquotients, 
where $K(X)$ is the quotient of $X$ by the involution $g_1$ (the Kummer 
variety). 

2. Let $X=\prod_{j=1}^NY_j$ be a product of generically twofold covers $Y_j$ 
of projective spaces (e.g., hyperelliptic curves) over $k$, at least one of 
which, for example $Y_1$, is a curve of genus $\le 1$. Then there are 
embeddings of $G$-representations $Z^d(k(X)\otimes_kF)\hookrightarrow
\bigoplus_{i=1}^NZ^d(k(\prod_{1\le j\le N,~j\neq i}Y_j)({\mathbb P}^{d_i})
\otimes_kF)^{1+\delta_{1i}}\hookrightarrow Z^d(k({\mathbb P}^d)\otimes_k
F)^{N+1}$, where $d_i=\dim Y_i$ and $d=\dim X$. In particular, 
${\rm JH}(X)={\rm JH}({\mathbb P}^d_k)$. 

\section{Differential forms and cohomologies}
\label{differential-forms} 
For any field extension $L|k$ let $H^q_{{\rm dR}/k,c}(L)$ be the 
image in $H^q_{{\rm dR}/k}(L):={\rm coker}[\Omega^{q-1}_{L|k}
\stackrel{d}{\longrightarrow}\Omega^q_{L|k,{\rm closed}}]$ of 
$\lim\limits_{\longrightarrow}H^q_{{\rm dR}/k}(X)$, where $X$ runs 
over smooth proper models of subfields in $L$ of finite type over $k$. 
If $L$ is algebraically closed, this is an admissible representation 
of $G_{L|k}$ over $k$, cf. the footnote on p.\pageref{opred-dopust}. This is 
a remarkable object and its structure reflects several, so far conjectural, 
relations between the cohomology and algebraic cycles. E.g., when 
$k={\mathbb C}$ the Hodge general conjecture is equivalent to the vanishing 
of any Hodge substructure in $H^q_{{\rm sing},c}(L)$ with $h^{q,0}=0$. In 
other words, it is equivalent to the vanishing of the subrepresentation 
${\mathfrak N}^1H^q_{{\rm sing},c}(F)$, which is by definition the maximal 
Hodge substructure in $H^q_{{\rm sing},c}(F)$ with $h^{q,0}=0$. 

The Hodge filtration on $\Omega^{\bullet}_{X|k}$ induces a descending 
filtration on $H^q_{{\rm dR}/k,c}(L)$ with the graded quotients 
$H^{p,q-p}_{L|k}=\lim\limits_{\longrightarrow}{\rm coker}[H^{p-1}
(D,\Omega^{q-p-1}_{D|k})\longrightarrow H^p(X,\Omega^{q-p}_{X|k})]$, 
where $(X,D)$ runs over the pairs consisting of a smooth proper 
variety $X$ with $k(X)\subset L$ and a normal crossing divisor $D$ 
on $X$ with smooth irreducible components. More particularly, 
$H^{q,0}_{L|k}=\Omega^q_{L|k,{\rm reg}}\subset H^q_{{\rm dR}/k,c}(L)$. 

It is easy to see that the restrictions of $\Omega^q_{F|k,{\rm reg}}$ and 
of $H^q_{{\rm sing},c}(F)/{\mathfrak N}^1$ to any compact subgroup $U$ of 
$G$ contain the same irreducible representations of $U$. Indeed, any finite 
group $\Gamma$ acting on a smooth proper complex variety $X$ induces 
automorphisms of the Hodge structure $H^q_{{\rm sing}}(X)$. For any 
irreducible representation $\rho$ of $\Gamma$ the element $p_{\rho}=
\frac{\deg\rho}{|\Gamma|}\sum_{\gamma\in\Gamma}{\rm tr}\rho(\gamma^{-1})
\gamma\in{\mathbb Q}[\Gamma]$ is a projector to the $\rho$-isotypical part. 
Then $p_{\rho}H^q_{{\rm sing}}(X)\subseteq H^q_{{\rm sing}}(X)_{{\mathbb Q}}$ 
is a Hodge substructure, and therefore, $\Gamma(X,\Omega^q_X)$ and 
$H^q_{{\rm sing}}(X)/{\mathfrak N}^1$ contain the same irreducible 
representations of $\Gamma$.

\begin{prop} Suppose that the cardinality of $k$ is at 
most continuum. Fix an embedding $\iota:k\hookrightarrow{\mathbb C}$ 
into the field of complex numbers. Then \begin{itemize} 
\item there is a non-canonical ${\mathbb Q}$-linear 
isomorphism $H^{p,q}_{F|k}\cong H^{q,p}_{F|k}$, and a 
${\mathbb C}$-anti-linear canonical isomorphism 
(depending on $\iota$) $H^{p,q}_{F|k}\otimes_{k,\iota}{\mathbb C}
\cong H^{q,p}_{F|k}\otimes_{k,\iota}{\mathbb C}$; 
\item the representation $H^n_{{\rm dR}/k,c}(L)$ of $G_{L|k}$ is 
semi-simple for any algebraically closed $L$ of a finite transcendence 
degree $n$ over $k$ (and in particular, $\Omega^n_{L|k,\text{{\rm reg}}}$ 
is also semi-simple). \end{itemize} \end{prop}
{\it Proof.} \begin{itemize}\item The complexification of the 
projection $F^pH^{p+q}_{{\rm dR}/k}(X)\longrightarrow\hspace{-3mm}\to 
H^q(X,\Omega^p_{X|k})$ identifies the space $F^pH^{p+q}_{{\rm dR}/k}(X)
\otimes_{k,\iota}{\mathbb C}\cap\overline{F^qH^{p+q}_{{\rm dR}/k}(X)
\otimes_{k,\iota}{\mathbb C}}$ with $H^q(X,\Omega^p_{X|k})
\otimes_{k,\iota}{\mathbb C}$, where $F^{\bullet}$ is the Hodge 
filtration. Then the complex conjugation on $H^{p+q}(X_{\iota}
({\mathbb C}),{\mathbb C})=H^{p+q}(X_{\iota}({\mathbb C}),{\mathbb R})
\otimes_{{\mathbb R}}{\mathbb C}$ identifies $H^q(X,\Omega^p_{X|k})
\otimes_{k,\iota}{\mathbb C}$ with $H^p(X_{\iota}({\mathbb C}),
\Omega^q_{X_{\iota}({\mathbb C})})=H^p(X,\Omega^q_{X|k})
\otimes_{k,\iota}{\mathbb C}$. 
\item The semi-simplicity of $H^n_{{\rm dR}/k,c}(L)$ is equivalent to 
the semi-simplicity of the representation ${\mathbb C}\otimes_{k,\iota}
H^n_{{\rm dR}/k,c}(L)=\bigoplus_{p+q=n}{\mathbb C}\otimes_{k,\iota}
H^{p,q}_{L|k}$ of $G_{L|k}$. For the latter note that there is a positive 
definite $G_{L|k}$-equivariant hermitian form $({\mathbb C}\otimes_{k,\iota}
H^{p,q}_{L|k})\otimes_{id,{\mathbb C},c}({\mathbb C}\otimes_{k,\iota}
H^{p,q}_{L|k})\longrightarrow{\mathbb C}(\chi)$, where $c$ is 
the complex conjugation and $\chi$ is the modulus of $G_{L|k}$, given 
by $(\omega,\eta)=\int_{X_{\iota}({\mathbb C})}i^{n^2+2q}\omega
\wedge\overline{\eta}\cdot[G_{L|k(X)}]$ for any $\omega,\eta\in 
H^{p,q}_{{\rm prim}}(X_{\iota}({\mathbb C}))={\mathbb C}\otimes
_{k,\iota}H^q_{{\rm prim}}(X,\Omega^p_{X|k})\subset{\mathbb C}
\otimes_{k,\iota}H^{p,q}_{L|k}$. Here $H^{p,q}_{{\rm prim}}
(X_{\iota}({\mathbb C}))$ denotes the subspace orthogonal to the 
sum of the images of all Gysin maps $H^{p-1,q-1}(D)\longrightarrow 
H^{p,q}(X_{\iota}({\mathbb C}))$ for all desingularizations $D$ 
of all divisors on $X_{\iota}({\mathbb C})$. \qed \end{itemize} 

\vspace{4mm}

\noindent
{\sl Acknowledgement.} {\small Both authors were supported in a part by 
INTAS (grant no. 05-1000008-8118). M.R. gratefully acknowledges the support 
by the Max-Planck-Institut f\"{u}r Mathematik in Bonn (Fall 2004), Alexander 
von Humboldt-Stiftung (2005 -- 2007), Pierre Deligne 2004 Balzan prize 
in mathematics (2006--2008), and RFBR (grants 
06-01-72550-CNRS-L{\fontencoding{OT2}\selectfont \_a}, 
07-01-92211-{\fontencoding{OT2}\selectfont NCNIL\_a}). 
Also he is grateful to Regensburg University for the hospitality and to 
Alexander von Humboldt-Stiftung for making his stay in Regensburg possible. 

}


\begin{thebibliography}{}
\bibitem[BZ]{bz} I.N.Bernstein, A.V.Zelevinsky, {\em Representations 
of the group $GL(n,F),$ where $F$ is a local non-Archimedean field.} 
Uspehi Mat. Nauk \textbf{31} (1976), no. 3(189), 5--70. 
\bibitem[Joh]{topos} P.T.Johnstone, {\bf Topos theory.} 
Academic Press, London, New York, San Francisco, 1977. 
\bibitem[Mi]{milne} J.S.Milne, \'{E}tale cohomology. Princeton 
Math.Series, 33. Princeton Univ.Press, Princeton, N.J., 1980. 
\bibitem[R1]{repr} M.Rovinsky, {\it Motives and admissible
representations of automorphism groups of fields.} Math. Zeit.,
{\bf 249} (2005), no. 1, 163--221, {\tt math.RT/0101170}.
\bibitem[R2]{pgl} M.Rovinsky, {\it Semilinear representations 
of PGL,} Selecta Math., {\bf 11} (2005), 491--522, {\tt math.RT/0306333}. 
\bibitem[R3]{max} M.Rovinsky, {\it On maximal proper 
subgroups of field automorphism groups.} {\tt math.RT/0601028}. 
\bibitem[R4]{rms} M.Rovinsky, 
{\it Automorphism groups of fields, and their representations,} 
Russian Math. Surveys, {\bf 62}:6 (2007), 1121--1186. 
\bibitem[SGA 4 I]{sga-4-1} {\bf Th\'{e}orie des topos et cohomologie 
\'{e}tale des sch\'{e}mas.} Tome 1 LNM 269, Springer-Verlag 1972. 
\bibitem[SGA 4 II]{sga-4-2} {\bf Th\'{e}orie des topos et cohomologie 
\'{e}tale des sch\'{e}mas.} Tome 2 LNM 270, Springer-Verlag 1972. 
\end{thebibliography}
\end{document}